\documentclass[hidelinks,onefignum,onetabnum]{siamart250211}

\ifpdf
  \DeclareGraphicsExtensions{.eps,.pdf,.png,.jpg}
\else
  \DeclareGraphicsExtensions{.eps}
\fi

\newsiamremark{remark}{Remark}
\newsiamremark{hypothesis}{Hypothesis}
\crefname{hypothesis}{Hypothesis}{Hypotheses}
\newsiamthm{claim}{Claim}
\newsiamremark{fact}{Fact}
\crefname{fact}{Fact}{Facts}

\usepackage{xcite}
\usepackage{xr}
\makeatletter
\newcommand*{\addFileDependency}[1]{
  \typeout{(#1)}
  \@addtofilelist{#1}
  \IfFileExists{#1}{}{\typeout{No file #1.}}
}
\makeatother

\usepackage{amsopn}

\usepackage{graphicx} 
\usepackage{srcltx,epsfig}
\usepackage{amsmath}
\usepackage{color}
\usepackage{float}
\usepackage{xcolor}
\usepackage{multirow}
\usepackage{verbatim}
\usepackage{pgfplots}
\pgfplotsset{compat=1.18} 
\usepackage{array}
\usepackage{booktabs} 
\usepackage{placeins}

\newcommand{\mF}{\mathcal{F}}
\usepackage{bm}
\usepackage{tikz}
\usetikzlibrary{shapes.geometric, arrows}
\usetikzlibrary{matrix}

\graphicspath{{images/}}

\usepackage{enumitem}
\usepackage{lipsum}
\usepackage{amsfonts}
\usepackage{epstopdf}
\usepackage{algorithmicx}
\usepackage{algorithm}
\usepackage{algpseudocode}
\renewcommand{\algorithmiccomment}[1]{\bgroup\hfill//~#1\egroup}

\algnewcommand\algorithmicpara{\textbf{Parameters:}}
\algnewcommand\Para{\item[\algorithmicpara]}
\algnewcommand{\AND}{\textbf{and} }
\algnewcommand{\OR}{\textbf{or} }
\algnewcommand{\Break}{\textbf{Break}}
\usepackage{amssymb,mathtools,mathrsfs}
\allowdisplaybreaks[2]
\usepackage{subfigure}
\usepackage{braket}
\usepackage{enumerate}

\newcommand{\ee}{{\rm e}}
\newcommand{\bx}{{\bm x}}
\newcommand{\mT}{{\mathcal{T}}}
\newcommand{\mFP}{{\mathcal{F}^{(p)}}}
\newcommand{\bzeta}{{\bm{\zeta}}}
\newcommand{\ii}{\mathrm{i}}
\newcommand{\dd}{\mathrm{d}}
\newcommand{\pp}{\partial}
\newcommand{\Span}{{\rm Span}}
\newcommand{\ra}{\rightarrow}

\DeclareMathOperator*{\argmin}{arg\,min}

\headers{Adaptive Hermite spectral method for Boltzmann equation}{S. Shao, Y. Wang, and J. Wu}

\title{An adaptive Hermite spectral method for the Boltzmann equation\thanks{\today.\funding{The first author was partially supported by the National Natural Science Foundation of China (Nos. 12325112, 12288101). The second author was partially supported by the Foundation of the President of China Academy of Engineering Physics (YZJJZQ2022017), the Science Challenge program (NO. TZ2025016), and the National Natural Science Foundation of China (Grant No. 12171026, and 12031013). This work was supported by the High-performance Computing Platform of Peking University.}}}

\author{
    Sihong Shao\thanks{CAPT, LMAM, and School of Mathematical Sciences, Peking University, Beijing, 100871, China (sihong@pku.edu.cn).} \and 
    Yanli Wang\thanks{Beijing Computational Science Research Center, Beijing, 100193, China (ylwang@csrc.ac.cn).} \and 
    Jie Wu\thanks{Center for Data Science, Peking University, Beijing, 100871, China (wujie5@stu.pku.edu.cn).}}

\ifpdf
\hypersetup{
  pdftitle={Adaptive Hermite spectral method for Boltzmann equation},
  pdfauthor={S. Shao, Y. Wang, and J. Wu}
}
\fi

\begin{document}

\maketitle

\begin{abstract}
We propose an adaptive Hermite spectral method for the three-dimensional velocity space of the Boltzmann equation guided by a newly developed frequency indicator. For the homogeneous problem, the indicator is defined by the contribution of high-order coefficients in the spectral expansion. For the non-homogeneous problem, a Fourier-Hermite scheme is employed, with the corresponding frequency indicator formulated based on distributions across the entire spatial domain. The adaptive Hermite method includes scaling and $p$-adaptive techniques to dynamically adjust the scaling factor and expansion order according to the indicator. Numerical experiments cover both homogeneous and non-homogeneous problems in up to three spatial dimensions. Results demonstrate that the scaling adaptive method substantially reduces $L^2$ errors at negligible computational cost, and the $p$-adaptive method achieves time savings of up to $74\%$.

\end{abstract}



\section{Introduction}
The Boltzmann equation, as a fundamental governing equation in kinetic theory, has wide-ranging applications in rarefied gas dynamics \cite{Bird1994, Dimarco2018} as well as in related fields such as semiconductor transport \cite{BenAbdallah1996}. However, its inherent high dimensionality and quadratic collision model pose significant challenges for numerical simulation. Historically, a variety of numerical schemes have been proposed, including both stochastic and deterministic approaches. Among the stochastic methods, the most representative one is the direct simulation Monte Carlo method \cite{Bird1994}, which has low complexity but usually suffers from stochastic fluctuations. The deterministic methods \cite{Dimarco2014}, including the discrete velocity models \cite{Goldstein1989}, Fourier spectral methods \cite{Pareschi1996, Pareschi2000b}, Hermite spectral method \cite{Wang2019, Hu2020}, mapped Chebyshev spectral method \cite{Hu2022a}, and so on, often face high complexity in high-dimensional simulations. Existing deterministic results for three-dimensional problems remain limited \cite{Dimarco2018}. In this work, we develop a scale-$p$-adaptive Hermite method for the Boltzmann equation that reduces error as well as computational cost, thereby enabling the efficient simulation of fully three-dimensional problems.

Adaptive scaling techniques for the Hermite spectral method have been studied for a long time. A strategy for selecting the scaling factor based on the support of the target function is proposed in \cite{Tang1993}, which inspires the subsequent works \cite{Shen2009, Hu2024a}. Recently, an adaptive spectral method, incorporating moving, scaling, and $p$-adaptive techniques, has been proposed for problems on unbounded domains in \cite{Xia2021, Xia2021a, Chou2023}. Specifically, a frequency indicator measuring the proportion of the high-order coefficients in the spectral expansion has been introduced for the scaling and $p$-adaptive techniques \cite{Xia2021, Xia2021a}, with successful applications to various problems \cite{Deng2025}. This progress inspires our work.

When the Hermite spectral method is applied to the velocity space of the Boltzmann equation, the distribution function is approximated directly in the unbounded domain $\mathbb{R}^3$. The method is naturally conservative and has spectral accuracy. However, a major limitation lies in the high computational cost of evaluating the collision operator, especially when simulating the non-homogeneous problem with large Knudsen numbers. In such regimes, the distribution function is far from local equilibrium, requiring a large expansion order in the Hermite method \cite{Hu2020} and thus leading to substantial computational expense. 

Another issue of the Hermite spectral method is that its numerical performance strongly depends on the choice of the expansion center \cite{Cai2021}, particularly the scaling factor. Appropriate choices could reduce the required expansion order, whereas unreasonable choices may even lose convergence. In existing studies, the determination of the expansion center usually depends on prior knowledge or physical properties of the specific problem. For example, the local temperature is chosen as the scaling factor in \cite{Cai2014d}. The scaling factor in the Hermite method is fixed as a manually chosen constant \cite{Wang2019, Hu2020}. An alternative approach based on the highest moment is proposed in \cite{Cai2021}, where the steady-state shock structure problem with high Mach number is studied with convergence results reported. Additionally, the physically dependent scaling factor is also employed in \cite{Filbet2022} to enhance convergence when applying the Hermite spectral method to Vlasov equations. Meanwhile, discussion of the adaptive expansion order in kinetic equations remains insufficient. A $p$-adaptive Hermite spectral method is proposed in \cite{Vencels2015} for the Vlasov equation, where the Hermite expansion order is adjusted based on the value of the highest order coefficient.

Inspired by the recent works \cite{Xia2021, Xia2021a}, we propose an adaptive Hermite spectral method for the Boltzmann equation to adaptively adjust the expansion order and scaling factor. For homogeneous problems, a frequency indicator is defined based on the high-order Hermite expansion coefficients. The $p$- and scaling adaptive methods are brought up by adjusting the expansion order and scaling factor according to the frequency indicator. A projection algorithm is constructed to realize the transform between Hermite expansions with different scaling centers, which maintains the conservation of the method. Then, a Fourier-Hermite scheme is employed for the non-homogeneous problem. Corresponding frequency indicator and adaptive methods are also constructed by taking the distributions on all collocation points into account. The proposed framework selects scaling factors more effectively, thereby improving the efficiency of the Hermite spectral method.

In the numerical experiments, three spatially homogeneous problems are first considered. The results confirm that the scaling adaptive method accurately captures the evolution of the distribution while significantly reducing the $L^2$ error. 
A mixed-Gaussian problem in which the non-adaptive Hermite method fails to converge is also examined, whereas the proposed scale-$p$-adaptive method achieves spectral convergence and reduces computational time by approximately half. Subsequently, spatially non-homogeneous problems are investigated, including two one-dimensional problems, the two-dimensional Taylor-Green vortex problem, and the three-dimensional Hexa-Gaussian perturbation problem. The reduction of error and the saving in computational time demonstrate the high efficiency of this adaptive method. 

The rest of this paper is organized as follows. The basic properties of the Boltzmann equation, the discretization of the microscopic velocity space, and the approximation of the quadratic collision model are first introduced in Sec.~\ref{sec:bg}. The $p$- and scaling adaptive method for the homogeneous problem is presented in Sec.~\ref{sec:4Adap}, and the adaptive Hermite method for the non-homogeneous problem is proposed in Sec.~\ref{sec:5Inhom}. Several numerical experiments are displayed in Sec.~\ref{sec:6Exp} to validate this adaptive Hermite spectral method. This paper concludes with some future works and details in Sec.~\ref{sec:conclusion} and supplementary materials.

\section{Preliminaries}\label{sec:bg}
This section begins with a brief introduction to the dimensionless Boltzmann equation. Then, the discretization of the microscopic velocity space in the framework of the Hermite spectral method is presented, followed by the detailed approximation of the Boltzmann collision term.

\subsection{Boltzmann equation}
In kinetic theory, the state of a large number of particles is described by a distribution function $f(t,\bm{x},\bm{v})$ defined on a six-dimensional phase space, where $t$ represents time, $\bm{x}$ and $\bm{v}$ represent the spatial position and the microscopic velocity, respectively. The Boltzmann equation governs the time evolution of $f$ by
\begin{equation}
\label{eq:2Bol}
    \frac{\pp f}{\pp t}+\bm{v}\cdot\frac{\pp f}{\pp\bm{x}}=\frac{1}{\varepsilon}Q[f,f],\qquad t\in\mathbb{R}^{+},\bm{x}\in\Omega\subset\mathbb{R}^3,\bm{v}\in\mathbb{R}^3,  
\end{equation}
where $\varepsilon$ is the Knudsen number. The Boltzmann collision term $Q$ accounts for binary collisions between particles:
\begin{equation}
\label{eq:2Q}
    Q[f,f](\bm{v}) =\int_{\mathbb{R}^3}\int_{\mathbb{S}^2} [f(\bm{v}_{*}^{\prime})f(\bm{v}^{\prime})-f(\bm{v}_{*})f(\bm{v})]B(|\bm{g}|,\chi)\,\dd\bm{\omega}\dd\bm{v}_{*}.
\end{equation}
In this context, $\bm{g} = \bm{v} - \bm{v}_*$ is the relative velocity, and $\chi = \arccos\left(\bm g \cdot \bm \omega/|\bm g|\right)$ is the scattering angle. The unit vector $\bm{\omega}$ denotes the direction of the post-collisional relative velocity. The post-collisional velocities $\bm{v}'$ and $\bm{v}_{*}'$ can be expressed as
\begin{equation}
\label{eq:2v}
    \bm v'=(\bm v+\bm v_{*}+|\bm v-\bm v_{*}|\bm\omega)/2,\qquad
    \bm v_{*}'=(\bm v+\bm v_{*}-|\bm v-\bm v_{*}|\bm\omega)/2.
\end{equation}
These expressions are derived from the microscopic conservation laws of momentum and energy. The specific form of the collision kernel $B(|\bm{g}|,\chi)$ depends on the collision model under consideration. In this work, we will focus on the inverse power law (IPL) model \cite{Bird1994}, where the inter-particle force is inversely proportional to the $\eta$-th power of their distance. The kernel then takes the following form:
{\small
\begin{equation}
\label{eq:2IPL}
    B(|\bm{g}|,\chi)=C_{\rm IPL}|\bm{g}|^{\frac{\eta-5}{\eta-1}}\frac{W_{0}}{\sin\chi}\left|\frac{\dd W_0}{\dd\chi}\right|,\qquad\eta>3,
\end{equation}}
where $C_{\rm IPL}$ is a constant given by nondimensionalization \cite{Hu2020} and $W_0$ is a positive impact parameter independent of $|\bm{g}|$. Detailed relation between $W_0$ and $\chi$ can be found in literatures \cite{Bird1994, Hu2020}. Within the IPL model, cases with $\eta > 5$ are referred to as ``hard potentials'', while those with $\eta < 5$ are called ``soft potentials''. When $\eta = 5$, the kernel $B$ becomes independent of $|\bm{g}|$; this special case is known as the ``Maxwell molecules'' model.

The collision operator $Q$ preserves mass, momentum, and energy as
\begin{equation}
    \int_{\mathbb{R}^3}\phi(\bm{v})Q[f,f]\dd\bm{v}=0,\qquad \phi(\bm{v})=1,\bm{v},|\bm{v}|^2/2.
\end{equation}
Another important property of the collision operator $Q$ is Boltzmann's $H$-theorem, which states that
\begin{equation}
\label{eq:H_theorem}
    \int_{\mathbb{R}^3}Q[f,f]\ln f\dd\bm{v}\leqslant 0.
\end{equation}
The equality holds if and only if $f$ reaches the equilibrium named Maxwellian
\begin{equation}
\label{eq:2M}
    \mathcal{M}(\bm{v})=\frac{\rho}{(2\pi\theta)^{3/2}}\exp\left(-\frac{|\bm{v}-\bm{u}|^2}{2\theta}\right),
\end{equation}
with the dimensionless macroscopic density $\rho$, velocity $\bm{u}$ and temperature $\theta$
\begin{align}
\label{eq:macro}
    \rho(t,\bm{x})=\int_{\mathbb{R}^3}f\dd\bm{v},\quad \bm{u}(t,\bm{x})=\frac{1}{\rho}\int_{\mathbb{R}^3}\bm{v}f\dd\bm{v},\quad \theta(t,\bm{x})=\frac{1}{3\rho}\int_{\mathbb{R}^3}|\bm{v}-\bm{u}|^2 f\dd\bm{v}.
\end{align}
In addition, other macroscopic variables of interest can be derived from $f$. Let $\bm{v}=(v_1,v_2,v_3)$ and $\bm{u}=(u_1,u_2,u_3)$, we define the stress tensor and heat flux as
\begin{align}
    &\sigma_{ij}(t,\bm{x})=\int_{\mathbb{R}^3}\left((v_i-u_i)(v_j-u_j)-\frac{\delta_{ij}}{3}|\bm{v}-\bm{u}|^2\right)f\dd\bm{v}, \qquad i,j=1,2,3,\label{eq:shear}\\
    &q_{i}(t,\bm{x})=\frac{1}{2}\int_{\mathbb{R}^3}(v_i-u_i)|\bm{v}-\bm{u}|^2f\dd\bm{v}, \qquad i=1,2,3.\label{eq:heatflux}
\end{align}

\subsection{Discretization of the microscopic velocity space}
In this section, we introduce the discretization of the microscopic velocity space in the framework of the Hermite spectral method. As a generalization of the Hermite spectral method for the Boltzmann operator \cite{Wang2019}, the scaled asymmetrically-weighted Hermite functions are adopted as the basis. The one-dimensional form reads
\begin{equation}
\label{eq:basis}
    \mathcal{H}_{k}^{\beta}(v)=\frac{\sqrt{\beta}}{(2\pi)^{1/4}\sqrt{2^{k}k!}}H_{k}\left(\frac{\beta v}{\sqrt{2}}\right)\ee^{-\frac{\beta^2 v^{2}}{2}},\qquad k\geqslant0,
\end{equation}
where $\beta > 0$ is a scaling factor. Here, $H_k$ denotes Hermite polynomials given by Rodrigues' formula:
\begin{equation}
\label{eq:Hermite}
    H_k(v)=(-1)^{k}\ee^{v^2}\frac{\dd^k}{\dd v^k}\left[\ee^{-v^2}\right].
\end{equation}
It is evident that $H_k$ shares the same odd-even parity with $k$, so does $\mathcal{H}_{k}^{\beta}$. We also define the scaled Hermite polynomial by dropping the exponential factor from $\mathcal{H}_{k}^{\beta}$:
\begin{equation}
\label{eq:test_func}
    H_{k}^{\beta}(v)=\frac{\sqrt{\beta}}{(2\pi)^{1/4}\sqrt{2^{k}k!}}H_{k}\left(\frac{\beta v}{\sqrt{2}}\right),\qquad k\geqslant 0.
\end{equation}
The three-dimensional basis functions are constructed via the tensor product:
\begin{align}
\label{eq:3H}
    \mathcal{H}_{k}^{\beta}(\bm{v})=\mathcal{H}_{k_1}^{\beta}(v_1)\mathcal{H}_{k_2}^{\beta}(v_2)\mathcal{H}_{k_3}^{\beta}(v_3), \qquad 
    H_{k}^{\beta}(\bm{v})=H_{k_1}^{\beta}(v_1)H_{k_2}^{\beta}(v_2)H_{k_3}^{\beta}(v_3),
\end{align}
where $k = (k_1, k_2, k_3)$ is a multi-index in three dimensions. It also satisfies the orthogonality:
\begin{equation}
\label{eq:orth_phi}
    \int_{\mathbb{R}^3}H_{l}^{\beta}(\bm{v})\mathcal{H}_{k}^{\beta}(\bm{v})\dd \bm{v}=\delta_{k_{1}l_{1}}\delta_{k_{2}l_{2}}\delta_{k_{3}l_{3}}.
\end{equation}
where $\delta_{kl}$ is the Kronecker Delta symbol. This framework leads naturally to a Petrov-Galerkin method. The distribution function $f$ is approximated by the finite expansion
\begin{equation}
\label{eq:3fN}
    f(t,\bx, \bm{v}) \approx f_{N}^{\bzeta,\beta}(t,\bx, \bm{v}) =\sum_{|k|=0}^{N}\hat{f}_{k}^{\bzeta,\beta}(t,\bm{x})\mathcal{H}_{k}^{\beta}(\bm{v}-\bzeta),
\end{equation}
where $N$ denotes the expansion order, $\bzeta\in\mathbb{R}^3$ denotes the moving center, and $|k| = k_1 + k_2 + k_3$ \cite{Wang2019, Hu2020}. The number of terms is $(N+1)(N+2)(N+3)/6=\mathcal{O}(N^3)$. The coefficients are computed as
\begin{equation}
\label{eq:coe_f}
    \hat{f}_{k}^{\bm{\zeta},\beta} =\int_{\mathbb{R}^3}H_{k}^{\beta}(\bm{v}-\bm{\zeta})f(\bm{v})\dd \bm{v},\qquad |k|\geqslant 0.
\end{equation}
Macroscopic physical variables can be derived from the expansion coefficients of $f$ as 
\begin{align}
\label{eq:3rhou}
    \rho=\frac{(2\pi)^{\frac{3}{4}}\hat{f}_{{0}}^{\beta}}{\beta^{\frac{3}{2}}},\quad
    u_i=\frac{\hat{f}_{{e}_i}^{\beta}}{\beta\hat{f}_{{0}}^{\beta}},
    \quad \theta=\frac{1}{3}\sum_{i=1}^{3}\left[\frac{1}{\beta^2}\left(\frac{\sqrt{2}\hat{f}_{2{e}_i}^{\beta}}{\hat{f}_{{0}}^{\beta}}+1\right)-u_{i}^2\right],
\end{align}
where $e_i\in\mathbb{N}^3$ refers to the multi-index vector with $1$ in the $i$-th position and $0$ elsewhere. For example, $e_1=(1,0,0)$.

For the basis function parameters, $\bzeta\in\mathbb{R}^3$ represents a moving center, and $\beta\in\mathbb{R}^+$ is a scaling factor that significantly influences the convergence of the truncated series. The parameters $\bzeta$, $\beta$, and the expansion order $N$ play a critical role in the accuracy of the approximation to $f$ in \eqref{eq:3fN}. Classically, there are several strategies for choosing $\bzeta$ and $\beta$. One common approach sets them according to physical variables by $\bzeta = \bm u$, and $\beta = 1 / \sqrt{\theta}$ \cite{Wang2019}. Alternatively, they may be chosen as global constants based on the physical properties of the problem \cite{Hu2020, Filbet2022}. However, these choices may be inadequate for distribution functions with significant spatial or temporal variation, potentially resulting in poor or even divergent spectral approximations.

In this work, we propose an adaptive method that dynamically adjusts the expansion parameters $\beta$ and $N$ based on the evolution of $f$, aiming for improvements in approximation quality. This is achieved by monitoring frequency indicators, which will be introduced in the following sections. As for the moving center $\bzeta$, we follow the method proposed in \cite{Hu2020}, where a globally uniform value is adopted. The development of an adaptive strategy for $\bzeta$ is left for future work.

Before introducing the adaptive method, we first present the approximation of the IPL collision model within the Hermite spectral framework in the following section.

\subsection{Approximation of the collision operator} \label{sec:32Q}
Due to the complex form of the quadratic collision term \eqref{eq:2Q}, the computational cost of approximating it is extremely high. The evaluation of the  Hermite expansion coefficients of the collision operator under the IPL model is provided in \cite{Wang2019}, and its generalization to the adaptive Hermite spectral method is rather straightforward. For simplicity, only a sketch of the approximation of the collision operator is given in this section. Detailed derivations and discussions are left to the supplementary material \ref{Supp:3}.

Consider approximating the collision term using a truncated series:
\begin{align}
    &Q[f_{N}^{\bm{\zeta},\beta},f_{N}^{\bm{\zeta},\beta}](\bm{v})\approx Q_{N}^{\bm{\zeta},\beta}(\bm v) =\sum_{|k|=0}^{N}\hat{Q}_{k}^{\bm{\zeta},\beta}\mathcal{H}_{k}^{\beta}(\bm{v}-\bm{\zeta}), \\ 
    &\hat{Q}_{k}^{\bm{\zeta},\beta}=\int_{\mathbb{R}^3}H_{k}^{\beta}(\bm{v}-\bm{\zeta})Q[f_{N}^{\bm{\zeta},\beta},f_{N}^{\bm{\zeta},\beta}](\bm{v})\dd \bm{v}, \quad 0\leqslant |k|\leqslant N.
\end{align}
The coefficients can be calculated by the following combined model
\begin{equation}\label{eq:5Qhat}
    \hat{Q}_{k}^{\bm{\zeta},\beta}=\begin{cases}
        \displaystyle\sum_{|i|=0}^{N_0}\sum_{|j|=0}^{N_0} A^{i,j;\beta}_{k}\hat{f}_{i}^{\bm{\zeta},\beta}\hat{f}_{j}^{\bm{\zeta},\beta}, & 0\leqslant|k|\leqslant N_0,\\
        \nu_{N_0}\left(\hat{\mathcal{M}}_{k}^{\bm{\zeta},\beta}-\hat{f}_{k}^{\bm{\zeta},\beta}\right), & N_0<|k|\leqslant N.
    \end{cases}
\end{equation}
Here, $N_0$ is a given constant order for the quadratic collision operator. $\left\lbrace A^{i,j;\beta}_{k}\right\rbrace$ is a precomputed tensor. $\nu_{N_0}$ denotes the spectral radius of the discrete linearized collision operator \cite{Wang2019}. $\hat{\mathcal{M}}_{k}^{\bm{\zeta},\beta}$ denotes the expansion coefficients of the local Maxwellian. The complexity will be $\mathcal{O}(N^9)$ if $N\leqslant N_0$, and $\mathcal{O}(N_0^9+N^3)$ if $N>N_0$.

It is known that by appropriately scaling the Hermite basis functions, one can achieve equal or better accuracy with a smaller expansion order $N$. In this way, a significant saving in computation time can be realized when $N<N_0$. This is precisely the goal of the adaptive method to be introduced in the next section.

\section{Adaptive Hermite spectral method for the homogeneous problem}
\label{sec:4Adap}
In this section, the adaptive Hermite method based on the idea presented in \cite{Xia2021, Xia2021a} is studied for the homogeneous setting. The moving center $\bzeta$ is fixed as the macroscopic velocity in the homogeneous problem, and it will be omitted temporarily in this section. The Petrov-Galerkin method for the homogeneous Boltzmann equation yields the following ODE system
\begin{equation}\label{eq:3ODE}
    \frac{\dd \hat{f}_{k}^{\beta}}{\dd t} =\hat{Q}_{k}^{\beta},\qquad 0\leqslant|k|\leqslant N.
\end{equation}
which can be discretized by a fourth-order Runge-Kutta method.

To construct the indicator for both the scaling and $p$-adaptive methods, we adopt a similar approach from \cite{Xia2021a}, where a frequency indicator is defined to measure the proportion of high-order terms in the truncated expansion. In general, this indicator has been shown to be closely related to the lower bound of interpolation error. Following this thought, an indicator applicable to both $p$-adaptivity and scaling is constructed. In the framework of the Hermite discretization \eqref{eq:3fN}, the desired indicator should measure the contribution of high-order terms in the entire expansion. It is therefore defined as 
\begin{equation}\label{eq:4Ftot}
    \mathcal{F}[f_{N}^{\beta}](t) \coloneq \left(\frac{\sum\limits_{|k|=N-\bar{N}+1}^{N}\left(\hat{f}_{k}^{\beta}(t)\right)^2}{\sum\limits_{|k|=0}^{N}\left(\hat{f}_{k}^{\beta}(t)\right)^2}\right)^{1/2},
\end{equation}
with $\bar{N}=\max\left\lbrace\lfloor N/3\rfloor,2\right\rbrace$. That is, the contribution of the highest $\bar{N}$-th order expansion coefficients is monitored.

In the following, the scaling and $p$-adaptive algorithms are developed in sequence. These are applied after each time step to maintain the indicator at a prescribed level by tuning the scaling factor $\beta$ and expansion order $N$. The projection algorithm for scaling adjustment is also introduced.

\subsection{The \texorpdfstring{$p$}{}-adaptive algorithm}
In this section, the detailed $p$-adaptive algorithm is introduced using the indicator \eqref{eq:4Ftot}. Based on the spectral convergence behavior with respect to the expansion order $N$, an increasing value of $\mathcal{F}[f_{N}^{\beta}]$ indicates that the current order $N$ is insufficient and should be increased. In contrast, a small value of $\mathcal{F}[f_{N}^{\beta}]$ suggests that the expansion order can be decreased by discarding the highest-order expansion coefficients without significantly compromising accuracy, thus reducing computational complexity.

In particular, for the spatially homogeneous Boltzmann equation, the distribution function $f$ converges to the Maxwellian \eqref{eq:2M} as time evolves. Consequently, with a suitable scaling factor such as $\beta = 1/\sqrt{\theta}$, only a few terms are needed in the Hermite expansion \eqref{eq:3fN} to approximate $f$ accurately when it is near equilibrium. Therefore, we can expect the expansion order $N$ to decrease gradually under the $p$-adaptive algorithm, leading to a substantial reduction in computational cost.

{\small
\begin{algorithm}
\caption{The $p$-adaptive algorithm.}
\label{alg:ada-p}
\begin{algorithmic}[1]
    \Require Coefficients $\hat{f}_{k}^{\beta}$, scaling factor $\beta$, expansion order $N$, reference indicator $\mathcal{F}^{(p)}$;
    \Para $N_{\max}$, $N_{\min}$, $\Delta N$, $0\leqslant\eta_{0}^{(p)}$, $0<\eta_{l}^{(p)}\leqslant1\leqslant\eta_{h}^{(p)}$, $\mathcal{F}_{h,0},\mathcal{F}_{l,0}\geqslant 0$.
    \State $\texttt{flag}\gets0$;
    \State $\mathcal{F}_{h}^{(p)}\gets\max(\eta_{h}^{(p)} \mathcal{F}^{(p)}, \mathcal{F}_{h,0})$, $\mathcal{F}_{l}^{(p)}\gets\max(\eta_{l}^{(p)} \mathcal{F}^{(p)}, \mathcal{F}_{l,0})$;
    \While{$\mathcal{F}[f_{N}^{\beta}](t^n)>\mathcal{F}_{h}^{(p)}$} \label{ln:pi}
    \State $\texttt{flag}\gets1$;
    \If{$N+\Delta N> N_{\max}$}
    \State \Break;
    \EndIf
    \State $N\gets N+\Delta N$;\Comment{increase $N$}
    \State $\hat{f}_{k}^{\beta}\gets0$ for $|k| \in [N-\Delta N+1,N]$;
    \EndWhile
    \While{$\texttt{flag}=0$ \AND $\mathcal{F}[f_{N}^{\beta}](t^n)<\mathcal{F}_{l}^{(p)}$ \AND $N-\Delta N<N_{\min}$}
    \State $\widetilde{N}\gets N-\Delta N$;
    \If{$\mathcal{F}[f_{\widetilde{N}}^{\beta}](t^n)<\mathcal{F}_{l}^{(p)}$}\label{ln:pd}
    \State $\hat{f}_{k}^{\beta}\gets0$ for $|k|\in [N-\Delta N+1,N]$;
    \State $N\gets\widetilde{N}$;\Comment{decrease $N$}
    \EndIf
    \EndWhile
    \If{$\texttt{flag}=1$ \OR ($N$ is modified \AND $\mathcal{F}[f_{N}^{\beta}](t^n)\geqslant\eta_0^{(p)} \mathcal{F}^{(p)}$)}
    \State $\mathcal{F}^{(p)}\gets\mathcal{F}[f_{N}^{\beta}](t^n)$;
    \EndIf
\end{algorithmic}
\end{algorithm}
}
The detailed algorithm is presented in Alg.~\ref{alg:ada-p}. Several key components of the algorithm are summarized below:
\begin{itemize}[leftmargin=2em]
\item At the start of the simulation, a reference indicator value is initialized using the initial distribution:
\begin{equation}
    \mathcal{F}^{(p)}=\mathcal{F}[f_{N}^{\beta}](0).
\end{equation}
The goal of the $p$-adaptive algorithm is to maintain the indicator near the reference value.

\item In the algorithm, $\mF^{(p)}_h$ and $\mF^{(p)}_l$ are the two thresholds primarily determined by the reference value. At time $t^n$, if the indicator $\mathcal{F}[f_{N}^{\beta}](t^n)$ exceeds the upper threshold $\mF_h^{(p)}$, the expansion order $N$ will be increased by $\Delta N$ to reduce the indicator. The newly added coefficients are set to zero. Conversely, if the current indicator is below the lower threshold $\mF_l^{(p)}$, the expansion order $N$ will be considered for reduction by $\Delta N$. Here, $\Delta N$ denotes the adaptive step size. Additionally, whenever the expansion order $N$ is modified, the reference value $\mathcal{F}^{(p)}$ is considered updated to the current indicator. The procedure for determining $\mF^{(p)}_h$ and $\mF^{(p)}_l$, as well as for updating $\mFP$, will be described in detail later.

\item $N_{\max}$ and $N_{\min}$ denote the maximum and minimum order allowed.If $N$ reaches $N_{\max}$ while the indicator still exceeds the threshold $\mathcal{F}_h^{(p)}$, the reference indicator value $\mathcal{F}^{(p)}$ will also be updated to the current indicator. This adjustment accounts for the fact that the original reference value can no longer be maintained, and we can only expect to keep it at the new level. 

\item During the simulation, it is possible for $\mFP$ to become too small. For instance, when $N$ increases rapidly, the indicator may temporarily drop to a very small value or even zero, resulting in a loss of reference significance. To address this issue, the thresholds $\mF_h^{(p)}$ and $\mF_l^{(p)}$ are decided as
\begin{equation}
    \label{eq:thre}
    \mF_h^{(p)} = \max\left(\eta_{h}^{(p)}\mathcal{F}^{(p)}, \mathcal{F}_{h,0}\right), \qquad \mF_l^{(p)} = \max\left(\eta_{l}^{(p)}\mathcal{F}^{(p)}, \mathcal{F}_{l,0}\right), 
\end{equation}
where $\mF_{h,0} > 0$ and $\mF_{l,0} > 0$ are two constant parameters used to prevent the thresholds from becoming too small. Parameters $\eta_{h}^{(p)} > 1$ and $\eta_{l}^{(p)} < 1$ determine the threshold and the sensitivity of the algorithm. Note that the algorithm becomes more sensitive as $\eta_{l}^{(p)}$ and $\eta_{h}^{(p)}$ are set closer to $1$.

Moreover, when updating the reference indicator value $\mFP$, we must also avoid a sudden drop. To that end, $\mFP$ is updated as 
\begin{equation}
    \label{eq:update_F_new}
     \mFP \leftarrow \mF[f_N^{\beta}](t^n), \qquad {\rm if}~ \mF[f_N^{\beta}](t^n) \geqslant \eta_0^{(p)}\mFP,
\end{equation}
where $\eta_0^{(p)}$ is a small positive constant.

\end{itemize}

In the simulation, parameters should be selected based on the specific problem. Here, we provide a representative configuration. The adaptive step size is set to $\Delta N = 1$. The lower bound of the expansion order is $N_{\min} = 2$, which ensures the conservation of mass, momentum, and energy. The upper bound, $N_{\max} = 24$, is determined by computational resources and the scale of the problem. The three indicator ratios are set as $\eta_0^{(p)}=0.01$, $\eta_l^{(p)}=0.3$, and $\eta_h^{(p)}=1.5$. The constant threshold values are set as $\mF_{l,0}=10^{-9}$ and $\mF_{h,0}=10^{-6}$. These parameters collectively govern the update of the reference indicator and the determination of the thresholds.

\subsection{The scaling adaptive algorithm}\label{sec:4scale}
When the distribution function $f$ is close to Maxwellian, the optimal scaling factor typically approaches $1/\sqrt\theta$. However, in certain situations, such as the early stages of the homogeneous case or in non-homogeneous cases near the kinetic regime, $f$ can be far from equilibrium. In these cases, our adaptive scaling algorithm may yield a more appropriate scaling factor. Compared to previous works \cite{Wang2019}, where the scaling factor is fixed as a constant, the following adaptive algorithm can better capture the evolutionary behavior of the numerical solution by adaptively adjusting the scaling factor.

Before presenting the adaptive scaling algorithm, we introduce a conservative projection operator $\mathcal{T}^{\beta\ra\beta'}$. This operator rescales $f_{N}^{\beta}$ from the current $\beta$ to a new $\beta'$. The relevant polynomial and function spaces are defined as follows:
\begin{equation}
\label{eq:space}
    P_N\coloneq{\rm Span}\{H_k^\beta(\bm{v}),0\leqslant|k|\leqslant N\},\qquad V_N^\beta\coloneq{\rm Span}\{\mathcal{H}_k^\beta(\bm{v}),0\leqslant|k|\leqslant N\}.
\end{equation}
Thus, the projection operator $\mathcal{T}^{\beta\ra\beta'}$, or the orthogonal projection from $V_{N}^{\beta}$ to $V_{N}^{\beta'}$, is defined as
\begin{equation}
\label{eq:transf}
    \mathcal{T}^{\beta\ra\beta'}[f_{N}^{\beta}](\bm{v})\coloneq \sum_{|l|=0}^{N}\hat{f}_{l}^{\beta'}\mathcal{H}_{l}^{\beta'}(\bm{v}),
\end{equation}
where the new coefficients are defined as
\begin{equation}\label{eq:4fp1}
    \hat{f}_{l}^{\beta'}\coloneq\int_{\mathbb{R}^3}H_{l}^{\beta^{\prime}}(\bm{v})f_{N}^{\beta}(\bm{v})\dd \bm{v}=\sum_{|k|=0}^{N}\hat{f}_{k}^{\beta}\int_{\mathbb{R}^3}H_{l}^{\beta^{\prime}}(\bm{v})\mathcal{H}_{k}^{\beta}(\bm{v})\dd \bm{v}.
\end{equation}
The detailed implementation of it will be presented later.

Similar to the previous $p$-adaptive algorithm, the scaling adaptive algorithm aims to maintain the indicator at a low level by adjusting the scaling factor. In the finite expansion \eqref{eq:3fN}, the scaling factor $\beta$ determines the approximating function space $V_{N}^{\beta}$, which directly affects the error. This algorithm also utilizes the same indicator \eqref{eq:4Ftot}, where a smaller value typically corresponds to a better approximation.

However, unlike in the $p$-adaptive algorithm, it is not straightforward to determine whether increasing or decreasing $\beta$ will reduce the indicator based solely on the current value of $\mathcal{F}[f_{N}^{\beta}]$. Therefore, a bidirectional line search method is employed to optimize $\beta$. The detailed procedure is presented in Alg.~\ref{alg:ada-s}, with several points listed below:
\begin{itemize}[leftmargin=2em]
\item The mesh of candidate $\beta$ values can be constructed in various ways. In this work, we restrict $\beta$ to an exponential mesh, defined using a constant parameter $q < 1$, which determines the resolution of the mesh:
\begin{equation}\label{eq:beta_mesh}
    S = \{\beta=q^m \mid m\in\mathbb{Z},~\beta_{\min}\leqslant q^m\leqslant \beta_{\max}\}.
\end{equation}
Here, $\beta_{\min}$ and $\beta_{\max}$ denote the lower and upper bounds of $\beta$, respectively.

\item The initial scaling factor $\beta$ is determined based on the initial distribution function $f(0,\bm{v})$, which is assumed to have an analytical form. We compute the adaptive Hermite expansion of $f(0,\bm{v})$ using all $\beta\in S$, and the value of $\beta$ with the smallest indicator is selected as the initial scaling factor:
\begin{equation}\label{eq:beta_0}
    \beta_0=\argmin_{\beta\in S}\mathcal{F}[f_{N}^{\beta}](0).
\end{equation}
The reference value is then initialized as
\begin{equation}
    \mathcal{F}^{(s)}=\mathcal{F}[f_{N}^{\beta_0}](0).
\end{equation}

\item In the algorithm, 
\begin{equation}
    \label{eq:thre_s}
    \mF_l^{(s)} = \eta_{l}^{(s)}\mathcal{F}^{(s)}, \qquad  \mF_h^{(s)} = \eta_{h}^{(s)}\mathcal{F}^{(s)}
\end{equation}
represent two thresholds determined by the reference value, where $\eta_l^{(s)} < 1$ and $\eta_h^{(s)} > 1$ are two parameters analogous to those in the $p$-adaptive algorithm. At time $t^n$, if the indicator $\mF[f_N^{\beta}](t^n)$ falls outside the reference range $[\mF_l^{(s)}, \mF_h^{(s)}]$, the algorithm attempts both decreasing and increasing the scaling factor $\beta$. The two new factors are chosen as 
\begin{equation}
    \label{eq:new_s}
    \beta_{\rm de} = \beta q, \qquad \beta_{\rm in} = \beta /q.
\end{equation}
Accordingly, the new coefficients are computed by transforming $f_N^{\beta}$ to the new scaling factors $\beta_{\rm de}$ and $\beta_{\rm in}$. Among the three scaling factors $\beta_{\rm de}, \beta$, and $\beta_{\rm in}$, the one that yields the minimum indicator value is selected as the new $\beta$, provided it does not reach the boundary of the mesh $S$:
\begin{equation}
    \label{eq:update_beta}
    \beta \leftarrow \argmin_{\beta'\in\{\beta_{\rm de}, \beta, \beta_{\rm in}\}}\mF[f_N^{\beta'}].
\end{equation}
The procedure defined by \eqref{eq:new_s} and \eqref{eq:update_beta} terminates when a local minimum of $\mF[f_N^{\beta}]$ is reached. Finally, if $\beta$ is modified during this procedure, the reference $\mathcal{F}^{(s)}$ will be updated 
\begin{equation}
    \label{eq:update_F_new_s}
     \mF^{(s)} \leftarrow \mF[f_N^{\beta}](t^n).
\end{equation}
\end{itemize}

To summarize the parameters, the lower and upper bounds of the scaling factor are $\beta_{\min}=0.15$ and $\beta_{\max}=2$, respectively. We set $q=0.9995$ to achieve a fine grid for $\beta$. The indicator ratios are chosen as $\eta_{l}^{(s)}=0.9995$, and $\eta_{h}^{(s)}=1.0005$. 
{\small 
\begin{algorithm}
\caption{The scaling adaptive algorithm.}
\label{alg:ada-s}
\begin{algorithmic}[1]
    \Require Coefficients $\hat{f}_{k}^{\beta}$, scaling factor $\beta$, expansion order $N$, reference indicator $\mathcal{F}^{(s)}$;
    \Para $\beta_{\max}$, $\beta_{\min}$, $0<\eta_{l}^{(s)}\leqslant1\leqslant\eta_{h}^{(s)}$, $0<q<1$.    
    \While{$\mathcal{F}[f_{N}^{\beta}](t^n)>\eta_{h}^{(s)} \mathcal{F}^{(s)}$ \OR $\mathcal{F}[f_{N}^{\beta}](t^n)<\eta_{l}^{(s)} \mathcal{F}^{(s)}$}
    \State $\beta_{\rm de} \gets \beta q$;
    \State $f_{N}^{\beta_{\rm de}} \gets \mT^{\beta\ra\beta_{\rm de}}[f_{N}^{\beta}]$;
    \State $\beta_{\rm in} \gets \beta/q$;
    \State $f_{N}^{\beta_{\rm in}} \gets \mT^{\beta\ra\beta_{\rm in}}[f_{N}^{\beta}]$;
    \If{$\mathcal{F}[f_{N}^{\beta_{\rm de}}](t^n)\leqslant\min\left\lbrace\mathcal{F}[f_{N}^{\beta}](t^n),\mathcal{F}[f_{N}^{\beta_{\rm in}}](t^n)\right\rbrace$ \AND $\beta_{\rm de}\geqslant\beta_{\min}$}
    \State $\beta \gets \beta_{\rm de}$;\Comment{decrease $\beta$}
    \State $f_{N}^{\beta}\gets f_{N}^{\beta_{\rm de}}$;
    \ElsIf{$\mathcal{F}[f_{N}^{\beta_{\rm in}}](t^n)\leqslant\min\left\lbrace\mathcal{F}[f_{N}^{\beta_{\rm de}}](t^n),\mathcal{F}[f_{N}^{\beta}](t^n)\right\rbrace$ \AND $\beta_{\rm in}\leqslant\beta_{\max}$}
    \State $\beta \gets \beta_{\rm in}$;\Comment{increase $\beta$}
    \State $f_{N}^{\beta}\gets f_{N}^{\beta_{\rm in}}$;
    \Else
    \State \Break;
    \EndIf
    \EndWhile
    \If{$\beta$ is modified}
    \State $\mathcal{F}^{(s)}\gets\mathcal{F}[f_{N}^{\beta}](t^n)$;
    \EndIf
\end{algorithmic}
\end{algorithm}
}

\paragraph{Implementation of the projection operator}
In this part, the implementation of the operator $\mT^{\beta\ra \beta'}$ in \eqref{eq:transf} is presented. Similar algorithms can be found in the literature \cite{Hu2020, Pagliantini2023, Issan2024}, though a different approach is proposed here. To compute the new coefficients, we define the integral
\begin{equation}\label{eq:4Tlk}
    T_{l,k}=\int_{\mathbb{R}}H_{l}^{\beta^{\prime}}(v)\mathcal{H}_{k}^{\beta}(v)\dd v,
\end{equation}
which possesses the properties stated in Prop.~\ref{prop:1}

\begin{proposition}\label{prop:1}
Let $T_{l,k}=0$ for $l<0$ or $k<0$, then it holds that 
\begin{align}
    T_{l,k}&=0, \qquad l < k\;{\rm or}\;2\nmid l+k, \label{eq:prop1}\\ 
    T_{l,k}&=\frac{\beta}{\beta'}\sqrt{\frac{k+1}{l+1}}T_{l+1,k+1}, \qquad 
    T_{l,l}=\left(\frac{\beta^{\prime}}{\beta}\right)^{l+\frac{1}{2}},\qquad k, l\geqslant 0, \label{eq:prop2}\\
    T_{l,k}&=\frac{\beta^\prime}{\beta}\sqrt{\frac{k+1}{l}}T_{l-1,k+1}+\frac{\beta^\prime}{\beta}\sqrt{\frac{k}{l}}T_{l-1,k-1}-\sqrt{\frac{l-1}{l}}T_{l-2,k}, \quad l\geqslant1,\; k\geqslant0, \label{eq:prop3}\\ 
    T_{l,k}&=\frac{\sqrt{l(k+1)}}{l-k}\left(\frac{\beta^\prime}{\beta}-\frac{\beta}{\beta^\prime}\right)T_{l-1,k+1},\qquad l \geqslant 1,\; k\geqslant 0,\; l\neq k.\label{eq:prop4}
\end{align}
\end{proposition}
The proof is provided in the supplementary material \ref{Supp:1}. $T_{l,k}$ can be calculated using Hermite-Gauss quadrature, which results in a computational complexity of $\mathcal{O}(N^3)$ for all $0\leqslant l, k\leqslant N$. However, thanks to the properties established above, an algorithm with reduced complexity $\mathcal{O}(N^2)$ can be devised for computing all $T_{l,k}$. Specifically, $T_{l,k}$ for $0\leqslant l,k \leqslant N$ forms an $(N+1)\times (N+1)$ lower triangular matrix with only about a quarter of the entries being nonzero. Therefore, we first evaluate the main diagonal entries $T_{l,l}$ for $0\leqslant l\leqslant N$, and then apply \eqref{eq:prop4} recursively to obtain the lower triangular entries.

With these $T_{l,k}$, the coefficients $\hat{f}_{k}^{\beta'}$ can be rewritten as 
\begin{equation}
    \hat{f}_{k}^{\beta'} =\sum_{|k|=0}^{N}\hat{f}_{k}^{\beta}T_{l_1,k_1}T_{l_2,k_2}T_{l_3,k_3},\label{eq:coef}
\end{equation}
which are computed through the following process:
\begin{equation}
\label{eq:cal_coe}
    \begin{aligned}        &g^{(1)}_{(l_1,k_2,k_3)}=\sum_{k_1=0}^{l_1}\hat{f}_{(k_1,k_2,k_3)}^{\beta}T_{l_1,k_1},\qquad 0\leqslant l_1+k_2+k_3\leqslant N,\\		&g^{(2)}_{(l_1,l_2,k_3)}=\sum_{k_2=0}^{l_2}g^{(1)}_{(l_1,k_2,k_3)}T_{l_2,k_2},\qquad 0\leqslant l_1+l_2+k_3\leqslant N,\\		&\hat{f}_{(l_1,l_2,l_3)}^{\beta^{\prime}}=\sum_{k_3=0}^{l_3}g^{(2)}_{(l_1,l_2,k_3)}T_{l_3,k_3},\qquad 0\leqslant l_1+l_2+l_3\leqslant N.
    \end{aligned}
\end{equation}
The computational complexity of \eqref{eq:cal_coe} is $\mathcal{O}(N^4)$. Thus, we have achieved an exact realization of the projection with overall complexity $\mathcal{O}(N^4)$. 

\begin{remark}
The projection operator $\mT^{\beta\rightarrow \beta'}$ satisfies that 
\begin{equation}\label{eq:4conser}
    \int_{\mathbb{R}^3}\phi(\bm{v})f_{N}^{\beta}(\bm{v})\dd \bm{v} =\int_{\mathbb{R}^3}\phi(\bm{v})\mathcal{T}^{\beta\ra\beta'}[f_{N}^{\beta}](\bm{v})\dd \bm{v}, \qquad\forall \phi\in P_{N}.
\end{equation}
This ensures the conservation of all macroscopic variables, including density, macroscopic velocity, and energy.
 
Although only the scaling transformation has been presented above, a generalized algorithm that simultaneously performs moving and scaling can be derived in a similar manner. In this case, the projection mapping a finite expansion centered at $[\bzeta, \beta]$ to a new one centered at $[\bzeta', \beta']$ can be calculated.
\end{remark}

So far, the adaptive algorithms for the homogeneous problem have been presented. The discussion of the adaptive algorithms for the non-homogeneous Boltzmann equation will be introduced in the next section.

\section{The non-homogeneous problem}\label{sec:5Inhom}
In this section, the adaptive Hermite spectral method is extended to the non-homogeneous Boltzmann equation. This work first considers the scenario with periodic boundary conditions, and extensions to other boundary conditions will be left for future work. For spatial discretization, the Fourier spectral method is employed, and the semi-discrete form of the non-homogeneous problem is presented in Sec.~\ref{sec:non_dis}. A new frequency indicator is introduced in Sec.~\ref{sec:non-indicator}. Finally, the complete adaptive scheme is summarized in Sec.~\ref{sec:4sch}.

\subsection{Semi-discrete scheme}\label{sec:non_dis}
This subsection introduces the discretization of both spatial and microscopic velocity spaces for the non-homogeneous Boltzmann equation, using the Fourier and adaptive Hermite spectral methods, respectively. For simplicity, the spatially one-dimensional case is presented. Generalization to three dimensions is straightforward. 

Let the spatial domain be $\Omega=[0, L]$. Under periodic boundary conditions, spatial discretization is performed using a Fourier basis. Assuming an even number $M$ as the expansion order in space, the basis functions are given by
\begin{equation}
\label{eq:fourier_basis}
    E_{l}(x)\coloneq\ee^{\ii lx\frac{2\pi}{L}},\qquad x\in\Omega, \qquad |l|\leqslant M/2.
\end{equation}
The associated collocation points are $x_j=j\frac{L}{M}$ for $j=0,1,\ldots,M-1$. Following the expansion in \eqref{eq:3fN}, the distribution function $f$ is approximated as
\begin{equation}\label{eq:5fNg}
    f(t,x,\bm v)\approx f_{N,M}^{\bzeta,\beta}(t,x,\bm v)
     =\sum_{|k|=0}^{N}\hat{f}_{k}^{\bzeta,\beta}(t, x)\mathcal{H}_{k}^{\beta}(\bm{v}-\bm{\zeta}),
\end{equation}
with coefficients defined by
\begin{equation}
\label{eq:coe_hf}
\begin{aligned}
    &\hat{f}_{k}^{\bzeta,\beta}(t, x)\coloneq\sum_{l=-M/2}^{M/2}\hat{g}_{k,l}^{\bzeta,\beta}(t)E_{l}(x), \\ 
    &\hat{g}_{k,l}^{\bzeta,\beta}(t) \coloneq\int_{\Omega}\int_{\mathbb{R}^3}E_l(x)H_{k}^{\beta}(\bm v-\bzeta) f(t,x,\bm v)\dd\bm v\dd x, \qquad  |l|\leqslant M/2.
\end{aligned}
\end{equation}
For the non-homogeneous problem, the local macroscopic variables can still be computed using \eqref{eq:3rhou}. Moreover, at the collision points $x_j, j = 0, \cdots, M-1$, the coefficients $\hat{f}_{k}^{\bzeta,\beta}(t, x)$ are approximated as 
\begin{equation}
    \label{eq:coll_f}
   \hat{f}_{k,j}^{\bzeta,\beta}(t) := \hat{f}_{k}^{\bzeta,\beta}(t,x_{j})=\sum_{l=-M/2}^{M/2}\hat{g}_{k,l}^{\bzeta,\beta}(t)E_{l}(x_j).
\end{equation}
In the adaptive method, the coefficients $\hat{f}_{k,j}^{\bzeta,\beta}(t)$, rather than $\hat{g}_{k,l}^{\bzeta,\beta}(t)$, are used to evolve the numerical solution due to their compatibility with the homogeneous formulation. Thanks to the Fast Fourier transform (FFT), the computational cost of converting between $\hat{f}_{k,j}^{\bzeta,\beta}(t)$ and $\hat{g}_{k,l}^{\bzeta,\beta}(t)$ is $\mathcal{O}(N^3 M\log M)$. 

With the expansion \eqref{eq:5fNg}, the convection term is approximated by 
\begin{equation}
\label{eq:5adv}
    v_1\frac{\pp f_{N,M}^{\bzeta,\beta}}{\pp x}(t,x_j,\bm v) = \sum_{|k|=0}^{N}\hat{h}_{k,j}^{\bzeta,\beta}(t)\mathcal{H}_{k}^{\beta}(\bm{v}-\bm{\zeta}),
\end{equation}
with
\begin{equation}
\label{eq:5h_hat}
    \hat{h}_{k,j}^{\bzeta,\beta} = \sum_{l=-M/2}^{M/2}\frac{2\pi\ii l}{L}\left(\frac{\sqrt{k_1}}{\beta}\hat{g}_{k-e_1,l}^{\bzeta,\beta}+\zeta_1 \hat{g}_{k,l}^{\bzeta,\beta}+\frac{\sqrt{k_1+1}}{\beta}\hat{g}_{k+e_d,l}^{\bzeta,\beta}\right)E_{l}(x_j).
\end{equation}
A detailed derivation is provided in the supplementary material \ref{Supp:2}. The total complexity of evaluating the convection term is $\mathcal{O}(N^3 M\log M)$. Additionally, the approximation of the collision operator follows the same approach described in Sec.~\ref{sec:32Q}, that is, \eqref{eq:5Qhat} is applied at each spatial collocation point $x_j$. Consequently, the semi-discrete form of the Boltzmann equation can be written as
\begin{equation}
\label{eq:semi-Boltzmann}
    \frac{\dd \hat{f}_{k,j}^{\bm{\zeta},\beta}(t)}{\dd t}+\hat{h}_{k,j}^{\bm{\zeta},\beta}(t) = \hat{Q}_{k,j}^{\bm{\zeta},\beta}(t),\quad 0\leqslant|k|\leqslant N,\; 0\leqslant j<M,
\end{equation}
where $\hat{h}_{k,j}^{\bm{\zeta},\beta}$ and $\hat{Q}_{k,j}^{\bm{\zeta},\beta}$ are given in \eqref{eq:5h_hat} and Sec.~\ref{sec:32Q}, respectively. In the simulation, a third-order Runge-Kutta method \cite{Shu1988} is utilized for time discretization. The time step length is chosen based on the CFL condition, which will be specified in Sec.~\ref{sec:6Exp}. 

This concludes the discretization of the Boltzmann equation using the Fourier-Hermite framework. The adaptive method is presented in the subsequent sections.

\subsection{Frequency indicator for non-homogeneous problem}
\label{sec:non-indicator}
Following the approach in Sec.~\ref{sec:4Adap}, a similar indicator to \eqref{eq:4Ftot} is constructed for the non-homogeneous case. With the spectral expansion \eqref{eq:5fNg} and the collocation points \eqref{eq:coll_f}, the frequency indicator is defined as 
\begin{equation}\label{eq:5ind}
    \mathcal{F}[f_{N,M}^{\bm{\zeta},\beta}](t) \coloneq \left(\frac{\displaystyle\sum_{|k|=N-\bar{N}+1}^{N}\sum_{j=0}^{M-1}\left(\hat{f}_{k,j}^{\bm{\zeta},\beta}(t)\right)^2}{\displaystyle\sum_{|k|=0}^{N}\sum_{j=0}^{M-1}\left(\hat{f}_{k,j}^{\bm{\zeta},\beta}(t)\right)^2}\right)^{1/2},
\end{equation}
where $\bar{N}=\max\left\lbrace\lfloor N/3\rfloor,2\right\rbrace$. This indicator, denoted by $\mathcal{F}[f_{N,M}^{\bzeta,\beta}]$, monitors the contribution of the highest-order $\bar{N}$ coefficients across the entire spatial expansion. In this way, the indicator defined by \eqref{eq:5ind} serves the same purpose as \eqref{eq:4Ftot}. 

In this work, the Fourier spectral method is utilized to discretize the spatial domain as a preliminary approach, enabling spectral accuracy in approximating the convection term. This choice avoids obscuring the effects of adaptive techniques in the microscopic velocity space. The indicator defined in \eqref{eq:5ind} can also be applied to other numerical schemes for the Boltzmann equation, such as the finite volume method, etc.

\begin{remark}
A more natural way to design the indicator involves the coefficients $\hat{g}_{k,l}^{\bzeta,\beta}(t)$, as in
{\small 
\begin{equation}\label{eq:5ind_rmk}
    \mathcal{F}[f_{N,M}^{\bm{\zeta},\beta}](t)\sim \left(\frac{\displaystyle\sum_{|k|=N-M+1}^{N}\sum_{l=-M/2}^{M/2}\left|\hat{g}_{k,l}^{\bzeta,\beta}(t)\right|^2}{\displaystyle\sum_{|k|=0}^{N}\sum_{l=-M/2}^{M/2}\left|\hat{g}_{k,l}^{\bzeta,\beta}(t)\right|^2}\right)^{1/2}.
\end{equation}}
Since the coefficients at spatial points are updated during the numerical simulation, $\hat{f}_{k,j}^{\bm{\zeta},\beta}(t)$ is used in \eqref{eq:5ind} instead of $\hat{g}_{k,l}^{\bzeta,\beta}(t)$. Moreover, by Parseval's formula 
{\small
\begin{equation}
\label{eq:parseval}
\sum_{l=-M/2}^{M/2}\left|\hat{g}_{k,l}^{\bm{\zeta},\beta}(t)\right|^2=\frac{1}{L}\int_{\Omega}\left(\hat{f}_{k}^{\bm{\zeta},\beta}(t,x)\right)^2\dd x \approx\frac{1}{M}\sum_{j=0}^{M-1}\left(\hat{f}_{k,j}^{\bm{\zeta},\beta}(t)\right)^2,
\end{equation}}
where the rectangular integration rule is utilized for the last approximation, which is exact for $E_{l}(x)$ with $|l|<M$ \cite{Shen2011}. Due to the rapid decay of Fourier coefficients, the quadrature error is negligible in practice, often reaching machine precision. As a result, \eqref{eq:5ind_rmk} is practically equivalent to \eqref{eq:5ind}.
\end{remark}

\subsection{Outline of the scheme}\label{sec:4sch}
The complete numerical scheme is presented in Alg.~\ref{alg:hom_scheme}, which incorporates both the $p$-adaptive algorithm and the scaling adaptive algorithm. These algorithms are similar to those used in the homogeneous case. The primary distinction lies in the choice of indicator. For the non-homogeneous problem, the adaptive algorithms are modified by replacing the indicator \eqref{eq:4Ftot} in Alg.~\ref{alg:ada-p} and \ref{alg:ada-s} with \eqref{eq:5ind}. The adjustment of the scaling factor $\beta$ is carried out by applying the projection $\mathcal{T}$ from Sec.~\ref{sec:4scale} to the local coefficients $\left\lbrace\hat{f}_{k,j}^{\bm{\zeta},\beta}\right\rbrace_k$ at all spatial collocation points $x_j$, $0\leqslant j<M$.

{\small
\begin{algorithm}
\caption{Scale-$p$-adaptive Hermite method for the Boltzmann equation.}
\label{alg:hom_scheme}
\begin{algorithmic}[1]
    \Require Initial distribution function $f(0,x,\bm{v})$, moving center $\bzeta\in\mathbb{R}^3$, velocity expansion orders $N_0,N$, spatial expansion order $M$, time step size $\Delta t$, collision tensor $\left\lbrace A^{i,j;1}_{k}\right\rbrace_{0\leqslant|i|,|j|,|k|\leqslant N_0}$;
    \State At the initial time step $t^0 = 0$, for a given initial distribution $f(0,x,\bm{v})$, moving center $\bzeta$, and expansion order $N$, determine the optimal scaling factor $\beta$ by \eqref{eq:beta_0} and compute the coefficients $\left(\hat{f}_{k,j}^{\bzeta,\beta}\right)^0$ for $0\leqslant|k|\leqslant N$, $0\leqslant j<M$. Initialize the reference indicators $\mathcal{F}^{(s)}$ and $\mathcal{F}^{(p)}$.
    
    \State \label{sch:42}At time step $t^n$, solve the ODE system \eqref{eq:semi-Boltzmann} using a third-order Runge-Kutta method (RK3) to obtain $\left(\hat{f}_{k,j}^{\bzeta,\beta}\right)^{n+1}$ for $0\leqslant|k|\leqslant N$, $0\leqslant j<M$.
    
    \State \label{sch:43}Check whether a scaling adjustment is required using Alg.~\ref{alg:ada-s}. If so, update the scaling factor $\beta$ and coefficients $\left(\hat{f}_{k,j}^{\bzeta,\beta}\right)^{n+1}$ for $0\leqslant|k|\leqslant N$, $0\leqslant j<M$.
    
    \State \label{sch:44}Check whether a $p$-adaptive adjustment is required using Alg.~\ref{alg:ada-p}. If so, update the expansion order $N$ and coefficients $\left(\hat{f}_{k,j}^{\bzeta,\beta}\right)^{n+1}$ for $0\leqslant|k|\leqslant N$, $0\leqslant j<M$. 
    
    \State If the final time is not reached, update $t^{n+1} =t^n+\Delta t$ and return to Step~\ref{sch:42}.
\end{algorithmic}
\end{algorithm}
} 
The scheme without Steps~\ref{sch:43},~\ref{sch:44} and fixed $\bzeta=\bm{u}$, $\beta=1/\sqrt{\theta}$ corresponds to the non-adaptive Hermite method described in \cite{Hu2020}. In the scale-$p$-adaptive method, both scaling and $p$-adaptivity are applied after each RK3 step to maintain an appropriate basis. The scaling factor $\beta$ is checked before the order $N$ to avoid unnecessary increases in $N$. Additionally, scaling-only or $p$-only adaptive schemes can be implemented by omitting Step~\ref{sch:44} or Step~\ref{sch:43}, respectively.

The computational complexity of each loop is dominated by the collision part in Step~\ref{sch:42}, which requires $\mathcal{O}(MN_0^9+MN^3)$ operations. For the adaptive components, the cost of evaluating the indicators is $\mathcal{O}(MN^3)$. The typical complexities of the scaling and $p$-adaptive steps are $\mathcal{O}(N_s MN^4)$ and $\mathcal{O}(N_p MN^3)$, respectively, where $N_s$ and $N_p$ denote the number of scaling and $p$-adjustments. These counts depend on the specific problem and parameter settings and are generally moderate. By dynamically adjusting the basis functions, the adaptive algorithm can improve the numerical accuracy at a relatively low cost.

\section{Numerical experiments}\label{sec:6Exp}
In this section, several numerical examples are tested to verify the effectiveness of the proposed adaptive techniques, including spatially homogeneous and non-homogeneous cases, where the dimension of the microscopic velocity space is fixed as three, while that of the spatial space varies from one to three. All tests are run on a single Intel Xeon Platinum 8358 processor. Multi-threading acceleration is utilized in the last two tests, as specified.

\subsection{Spatially homogeneous case}
\label{sec:0D}
This section presents three numerical examples of the spatially homogeneous Boltzmann equation solved using the adaptive scheme. The first example is taken from \cite{Wang2019}, while the second is a Quad-Gaussian problem. The third example is a mixed-Gaussian initial value problem that is intractable for the non-adaptive method. For simplicity, the dimensionless coefficient $\varepsilon$ in \eqref{eq:2Bol} and $C_{\rm IPL}$ in \eqref{eq:2IPL} are both set to 1, i.e., $C_{\rm IPL} = \varepsilon = 1$.

\subsubsection{BKW solution}
\label{sec:Eg01}
In this example, the performance of the scaling adaptive scheme is investigated with the BKW solution problem. The BKW solution has an analytical form \eqref{eq:BKW} for $\eta=5$ and has been studied in \cite{Bobylev1975, Krook1977, Wang2019}.
{\small 
\begin{equation}\label{eq:BKW}
\begin{aligned}
        f(t,\bm v)& =\frac{1}{(2\pi S(t))^{3/2}}\exp\left(-\frac{|\bm v|^2}{2S(t)}\right)\left(\frac{1-S(t)}{2(S(t))^2}|\bm v|^2+\frac{5S(t)-3}{2S(t)}\right), \\
         S(t)&=1-0.4\exp\left(\pi\tilde{B}_2^5t/3\right), \qquad \tilde{B}_2^5\approx-0.6543.
\end{aligned}
\end{equation}
}
Explicit definition of $\tilde{B}_2^5$ can be found in \cite{Wang2019}. As time $t$ increases, $S(t)$ approaches 1, and the distribution function $f$ converges to a Maxwellian. The temperature in this problem remains constant, $\theta = 1$. In the simulation, the time step is set to $\Delta t = 0.001$, and the parameters used in the scaling adaptive scheme are listed in Tab.~\ref{tab:01-1}. This example tests expansion orders $N = N_0 = 4, 6, 8$, where $N_0$ denotes the order of the quadratic collision term defined in \eqref{eq:5Qhat}. Additionally, the moving center $\bzeta$ is set to 0 since the distribution function $f$ is an even function.

\begin{table}
    \caption{(BKW solution in Sec.~\ref{sec:Eg01}) Parameters utilized in the scaling adaptive algorithm.}
    \label{tab:01-1}
    \centering
    \def\arraystretch{1.3}
    {\footnotesize
    \begin{tabular}{c|ccccc}
    parameter & $\beta_{\min}$ & $\beta_{\max}$ & $q$ & $\eta_{l}^{(s)}$ & $\eta_{h}^{(s)}$ \\ \hline
    usage & \eqref{eq:beta_mesh} & \eqref{eq:beta_mesh} & \eqref{eq:beta_mesh} & \eqref{eq:thre_s} & \eqref{eq:thre_s} \\
    value & 0.15 & 2.0 & 0.9995 & 0.9995 & 1.0005
    \end{tabular}
    }
\end{table}

At the beginning of the simulation, the initial scaling factor $\beta_0$ in \eqref{eq:beta_0} is selected by the scaling adaptive scheme. We can observe that if the scaling factor is chosen as $\beta_{\rm opt}(t) = 1/\sqrt{S(t)}$, the analytical solution \eqref{eq:BKW} can be represented exactly by the spectral expansion \eqref{eq:3fN} with only $N = 2$. Therefore, the scaling factor determined by the scaling adaptive scheme is compared with the optimal factor $\beta_{\rm opt}$. At $t = 0$, the initial $\beta_0$ closely matches the optimal value $\beta_{\rm opt}(0) \approx 1.2910$, demonstrating the accuracy of the adaptive scheme in selecting the initial scaling factor. The $L^2$ error between the numerical solution obtained using the scaling adaptive scheme and the analytical solution at $t = 0$ is reported in Tab.~\ref{tab:01-2}. For comparison, the non-adaptive Hermite method with a constant scaling factor $\beta$ is used as a reference. In this method, $\beta$ is set according to the temperature, i.e., $\beta \equiv 1/\sqrt{\theta}$. The corresponding $L^2$ error at $t = 0$ is also listed in Tab.~\ref{tab:01-2}. Results show that the initial $L^2$ error of the scaling adaptive scheme is significantly smaller than that of the non-adaptive method, validating the effectiveness of the adaptive approach.

\begin{figure}
    \centering
    \subfigure[$L^2$ error\label{fig:01-1a}]{
        \includegraphics[height=7.4\baselineskip]{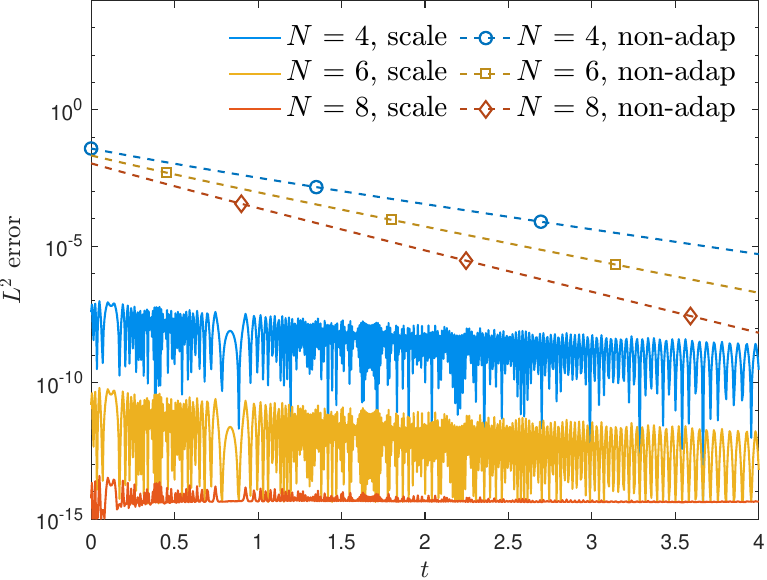}
    }\hspace{-0.7em}
    \subfigure[Scaling factor $\beta$\label{fig:01-1b}]{
        \includegraphics[height=7.4\baselineskip]{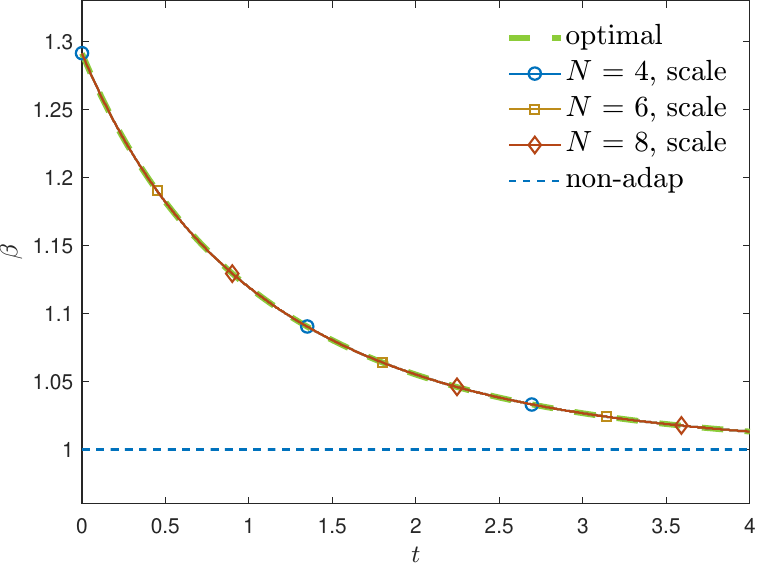}
    }\hspace{-0.7em}
    \subfigure[Frequency indicator $\mathcal{F}$ \label{fig:01-1c}]{
        \includegraphics[height=7.4\baselineskip]{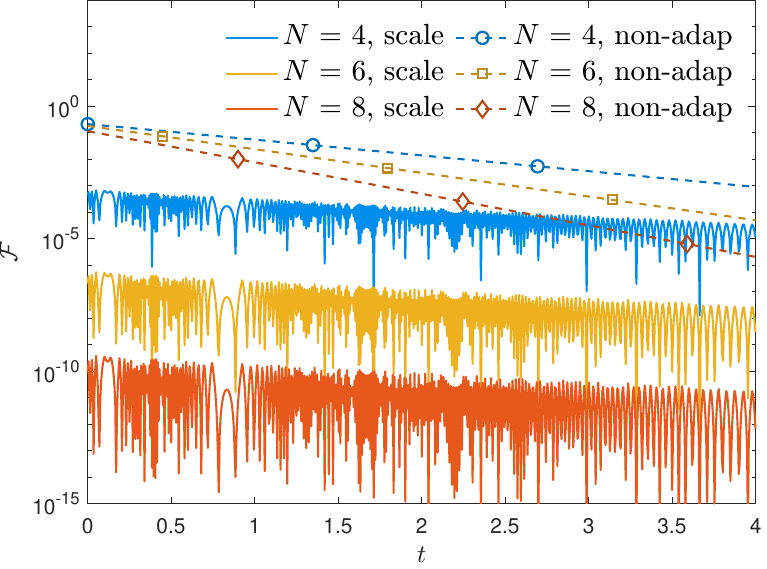}
    }
    
    \caption{(BKW solution in Sec.~\ref{sec:Eg01}) Comparison between the scaling adaptive and non-adaptive methods with different expansion order $N$. (a) Evolution of the $L^2$ error between the numerical and the analytical solutions. (b) Evolution of the scaling factor $\beta$. (c) Evolution of the frequency indicator $\mF$.}
    \label{fig:01-1}
\end{figure}

The evolution of the $L^2$ error, the numerical scaling factor $\beta$, and the indicator $\mF$ are shown in Fig.~\ref{fig:01-1}, respectively. As illustrated in Fig.~\ref{fig:01-1a}, for a given expansion order $N$, the $L^2$ error achieved by the scaling adaptive scheme is much lower than that of the non-adaptive method. Notably, with $N = 4$, the error reaches the order of $10^{-8}$. To compare two methods quantitatively, the $L^2$ error at $t = 4$ is listed in Tab.~\ref{tab:01-2}, further demonstrating the improved accuracy of the scaling adaptive scheme over the non-adaptive one. Fig.~\ref{fig:01-1b} shows that the dynamically adjusted scaling factor $\beta$ closely matches the optimal value $\beta_{\rm opt}(t)$, which indicates that the scaling adaptive scheme can track the evolution of the numerical solution accurately. Besides, the behavior of the indicator in Fig.~\ref{fig:01-1c} aligns with that of the $L^2$ error, which confirms the reliability and effectiveness of the proposed frequency indicator.

The average computational time of each time step for the scaling adaptive method is listed in Tab.~\ref{tab:01-3}, where the time for the non-adaptive method, and the cost of the adaptive adjustment are all illustrated. It reveals that the computational cost of the scaling adaptive method is nearly identical to that of the non-adaptive method. Moreover, the cost introduced by the adaptive adjustment is relatively small, highlighting the high efficiency of the proposed scaling adaptive method. 

\begin{table}[hptb]
    \caption{(BKW solution in Sec.~\ref{sec:Eg01}) $L^2$ error between the numerical and the analytical solutions for $N = 4$, $6$, and $8$ at $t = 0$ and $t = 4$. Numerical solutions are obtained using the scaling adaptive and non-adaptive methods.}
    \label{tab:01-2}
    \centering
    \def\arraystretch{1.3}
    {\footnotesize
    \begin{tabular}{c||ccc|ccc}
     & \multicolumn{3}{c|}{$t=0$} & \multicolumn{3}{c}{$t=4$} \\ 
     & $N=4$ & $N=6$ & $N=8$ & $N=4$ & $N=6$ & $N=8$ \\ \hline
    scale & 3.84E$-$08 & 1.55E$-$11 & 5.85E$-$15 & 2.36E$-$10 & 4.53E$-$14 & 4.44E$-$15 \\
    non-adap & 3.80E$-$02 & 2.12E$-$02 & 1.09E$-$02 & 5.12E$-$06 & 1.97E$-$07 & 6.76E$-$09
    \end{tabular}
    }
\end{table}

\begin{table}[hptb]
    \caption{(BKW solution in Sec.~\ref{sec:Eg01}) Average CPU time per time step (in seconds). Here, $T_{\rm non}$ and $T_{\rm adap}$ refer to the CPU time of the non-adaptive and scaling adaptive methods, respectively. $T_{\rm ind}$ refers to the CPU time spent on the adaptive adjustment step, which is a component of $T_{\rm adap}$.}
    \label{tab:01-3}
    \centering
    \def\arraystretch{1.3}
    {\footnotesize
    \begin{tabular}{c||cccc}
     & $T_{\rm non} $ & $T_{\rm adap}$ & $T_{\rm ind}$ & $T_{\rm ind}/T_{\rm adap}$ \\ \hline
    $N=4$ & 4.30E$-$5 & 4.57E$-$5 & 2.50E$-$6 & 5.47\% \\
    $N=6$ & 7.80E$-$4 & 7.89E$-$4 & 5.78E$-$6 & 0.73\% \\
    $N=8$ & 6.63E$-$3 & 6.71E$-$3 & 1.43E$-$5 & 0.21\%
    \end{tabular}
    }
\end{table}

\subsubsection{Quad-Gaussian problem}\label{sec:Eg03}
To further examine the effect of the scaling adaptive method, a symmetric four-component Gaussian mixture problem is considered in this example. The initial condition consists of the sum of four Maxwellian distributions with identical density and temperature but different velocities, given by
{\small 
\begin{gather}
    f(0,\bm{v})= \sum_{i = 1}^4 \frac{\rho_i}{(2\pi\theta_i)^{3/2}}\exp\left(-\frac{|\bm{v}-\bm{u}_i|^2}{2\theta_i}\right),\\
    \rho_i = 1/4,\qquad \theta_i = 1/3, \qquad i=1,\dots, 4,\\
    \bm{u}_1=-\bm{u}_2=(\sqrt{2},0,0),
    \qquad \bm{u}_3=-\bm{u}_4=(0,\sqrt{2},0).
\end{gather}
}
The temperature remains constant as $\theta\equiv1$. In the simulation, the time step length is set to $\Delta t=0.01$, and the parameters of the scaling adaptive scheme are the same as those in Tab.~\ref{tab:01-1}. Expansion orders $N = N_0 = 10, 12, 14$ are tested. The moving center $\bzeta$ is fixed at $0$. The reference solution is obtained by the scaling adaptive method with $N = N_0 = 20$.

\begin{figure}
    \centering
    \subfigure[$L^2$ error\label{fig:03-1a}]{
        \includegraphics[height=7.4\baselineskip]{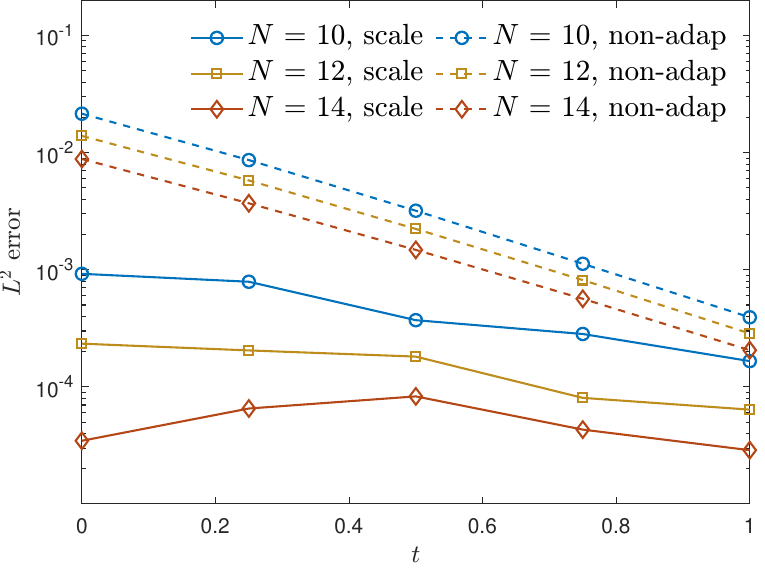}
    }\hspace{-0.7em}
    \subfigure[scaling factor $\beta$\label{fig:03-1b}]{
        \includegraphics[height=7.4\baselineskip]{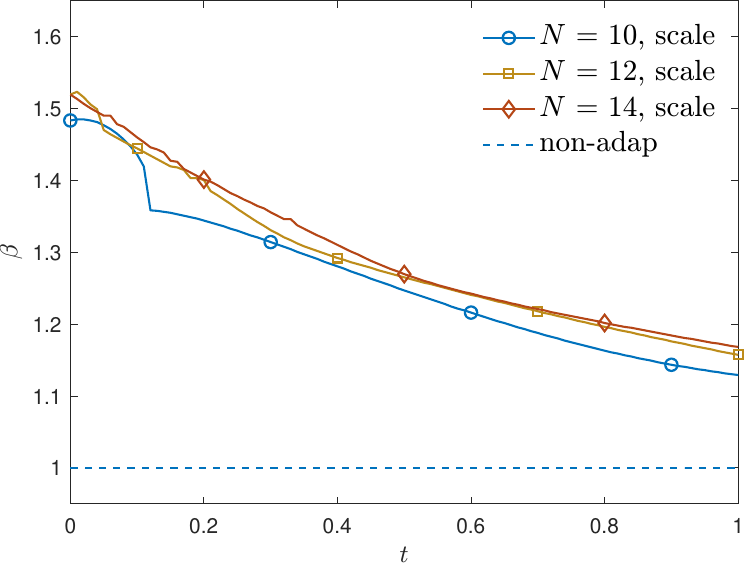}
    }\hspace{-0.7em}
    \subfigure[frequency indicator $\mathcal{F}$\label{fig:03-1c}]{
        \includegraphics[height=7.4\baselineskip]{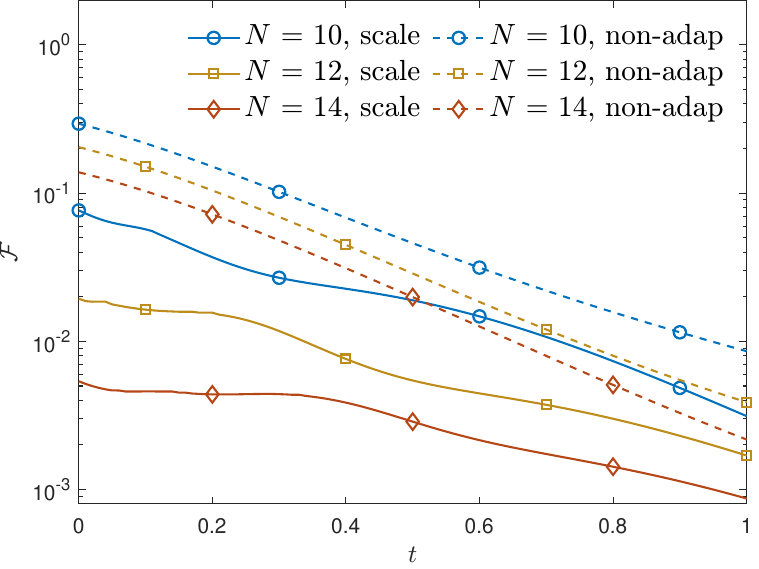}
    }
    
    \caption{(Quad-Gaussian problem in Sec.~\ref{sec:Eg03}) Comparison between the scaling adaptive and non-adaptive methods with different expansion order $N$. (a) Evolution of the $L^2$ error between the numerical and the reference solution. (b) Evolution of the scaling factor $\beta$. (c) Evolution of the frequency indicator $\mF$.}
    \label{fig:03-1}
\end{figure}

At the beginning, the initial scaling factor $\beta_0$ in \eqref{eq:beta_0} is determined by the scaling adaptive scheme. Consequently, the initial $L^2$ errors for $N = 10, 12$ and $14$, which are listed in Tab.~\ref{tab:03-1}, are substantially smaller than those of the non-adaptive method. The non-adaptive method here refers to the Hermite spectral method with the scaling factor fixed as the temperature $\beta \equiv 1 / \sqrt{\theta}= 1$. The evolution of the $L^2$ error is shown in Fig.~\ref{fig:03-1a}, and the corresponding values at $t = 1$ are also given in Tab.~\ref{tab:03-1}. It can be seen that the scaling adaptive method consistently yields significantly lower $L^2$ errors than the non-adaptive method. Notably, with $N = 10$, the scaling adaptive method achieves a smaller $L^2$ error than the non-adaptive method with $N = 14$. Furthermore, the time evolution of the scaling factor $\beta$ and the frequency indicator $\mF$ is displayed in Fig.~\ref{fig:03-1b} and Fig.~\ref{fig:03-1c}, respectively, showing similar trends between the error evolution and the frequency indicator.

\begin{figure}
    \centering
        \subfigure[$t=0$\label{fig:03-2a}]{
        \includegraphics[width=0.32\linewidth]{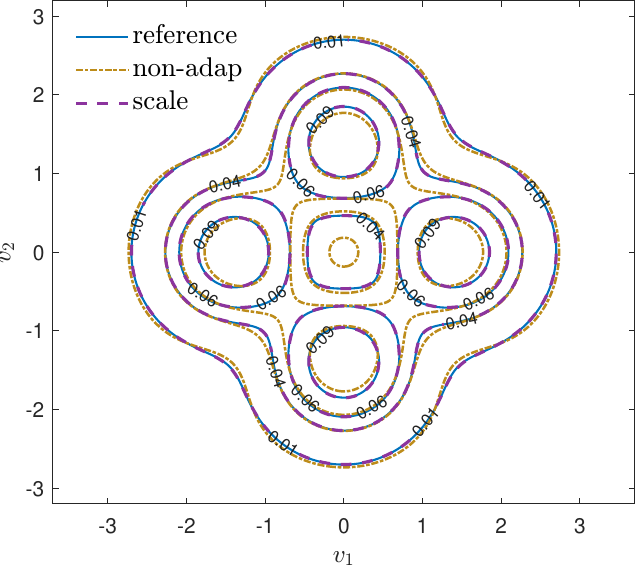}
    }\hspace{-0.7em}
    \subfigure[$t=0.25$\label{fig:03-2b}]{
        \includegraphics[width=0.32\linewidth]{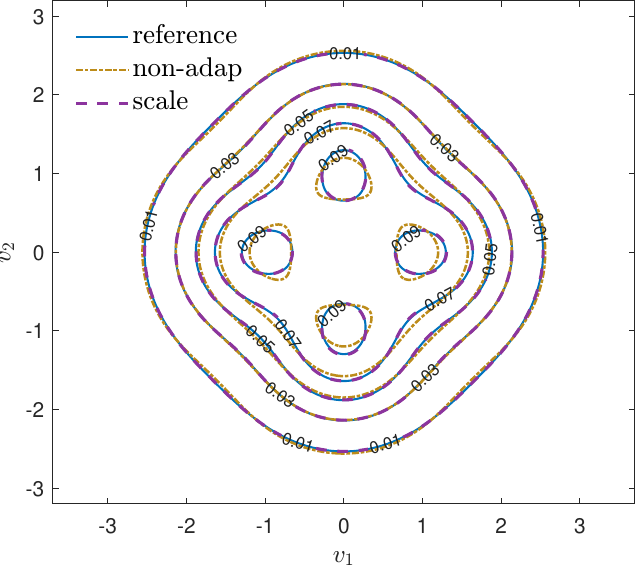}
    }\hspace{-0.7em}
    \subfigure[$t=1$\label{fig:03-2c}]{
        \includegraphics[width=0.32\linewidth]{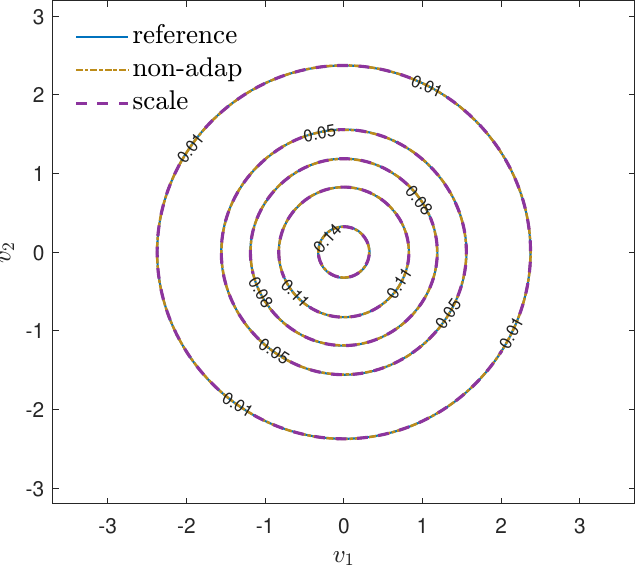}
    }
    
    \caption{(Quad-Gaussian problem in Sec.~\ref{sec:Eg03}) Marginal distribution functions for $N = 10$ at different times. Results from the scaling adaptive and non-adaptive methods are compared with the reference solution. (a) $t = 0$. (b) $t = 0.25$. (c) $t = 1$.}
    \label{fig:03-2}
\end{figure}

The marginal distribution functions for $N = 10$ at $t = 0$, $0.25$, and $1$ are illustrated in Fig.~\ref{fig:03-2}, along with the reference solution obtained by the scaling adaptive method with $N = 20$ and the numerical solution from the non-adaptive method with $N = 10$. At $t = 0$ and $0.25$, the numerical solution obtained by the scaling adaptive method agrees closely with the reference solution, whereas the non-adaptive method exhibits notable deviations. By $t = 1$, the distribution function approaches the Maxwellian, and both methods produce results that closely match the reference solution.

Further, the average computational time per time step for the scaling adaptive method is reported in Tab.~\ref{tab:03-2}, along with the corresponding times for the non-adaptive method and the cost of the adaptive adjustment. As in Sec.~\ref{sec:Eg01}, the computational costs of the scaling adaptive and non-adaptive methods are nearly identical, and the proportion attributable to the adaptive procedure is negligible.

\begin{table}[hptb]
    \caption{(Quad-Gaussian problem in Sec.~\ref{sec:Eg03}) $L^2$ error between the numerical and the reference solution for $N = 10$, $12$, and $14$ at $t = 0$ and $t = 1$. Numerical solutions are obtained using the scaling adaptive and non-adaptive methods.}
    \label{tab:03-1}
    \centering
    \def\arraystretch{1.3}
    {\footnotesize
    \begin{tabular}{c||ccc|ccc}
     & \multicolumn{3}{c|}{$t=0$} & \multicolumn{3}{c}{$t=1$} \\
     & $N=10$ & $N=12$ & $N=14$ & $N=10$ & $N=12$ & $N=14$ \\ \hline
    scale & 9.23E$-$4 & 2.34E$-$4 & 3.46E$-$5 & 1.66E$-$4 & 6.39E$-$5 & 2.88E$-$5 \\
    non-adap & 2.15E$-$2 & 1.39E$-$2 & 8.83E$-$3 & 3.93E$-$4 & 2.86E$-$4 & 2.06E$-$4
    \end{tabular}
    }
\end{table}

\begin{table}[hptb]
    \caption{(Quad-Gaussian problem in Sec.~\ref{sec:Eg03}) Average CPU time per time step (in seconds). Here, $T_{\rm non}$ and $T_{\rm adap}$ refer to the CPU time of the non-adaptive and scaling adaptive methods, respectively. $T_{\rm ind}$ refers to the CPU time spent on the adaptive adjustment step, which is a component of $T_{\rm adap}$.}
    \label{tab:03-2}
    \centering
    \def\arraystretch{1.3}
    {\footnotesize
    \begin{tabular}{c||cccc}
     & $T_{\rm non} $ & $T_{\rm adap}$ & $T_{\rm ind}$ & $T_{\rm ind}/T_{\rm adap}$ \\ \hline
    $N=10$ & 4.51E$-$2 & 4.52E$-$2 & 7.90E$-$5 & 0.17\% \\
    $N=12$ & 1.99E$-$1 & 2.03E$-$1 & 1.36E$-$4 & 0.07\% \\
    $N=14$ & 6.15E$-$1 & 6.15E$-$1 & 2.07E$-$4 & 0.03\%
    \end{tabular}
    }
\end{table}

\subsubsection{Mixed-Gaussian problem}
\label{sec:Eg04}
This example examines an asymmetric two-component Gaussian mixture problem, where the initial condition is given by
{\small
\begin{gather}
    f(0,\bm{v})=\frac{\rho_1}{(2\pi\theta_1)^{3/2}}\exp\left(-\frac{|\bm{v}-\bm{u}_1|^2}{2\theta_1}\right)+\frac{\rho_2}{(2\pi\theta_2)^{3/2}}\exp\left(-\frac{|\bm{v}-\bm{u}_2|^2}{2\theta_2}\right),\label{eq:5BiGa}\\
    \rho_1=1-\kappa, \qquad \bm{u}_1=\left(\sqrt{5/3}\mu,0,0\right), \qquad  \theta_1=1, \qquad \kappa = 0.05,\qquad \mu = 4, \\
    \rho_2=\kappa\frac{4\mu^2}{\mu^2+3},\qquad \bm{u}_2=\left(\sqrt{\frac{5}{3}}\frac{\mu^2+3}{4\mu},0,0\right), \qquad \theta_2=\frac{(5\mu^2-1)(\mu^2+3)}{16\mu^2}.
\end{gather}
}
Such a condition is discussed in \cite{Cai2021}. For this mixed-Gaussian problem, the macroscopic velocity and temperature are
{\small
\begin{equation*}
    \bm u = \frac{\rho_1\bm u_1+\rho_2\bm u_2}{\rho_1+\rho_2}\approx(4.62,0,0),\quad
    \theta = \frac{\rho_1\theta_1+\rho_2\theta_2}{\rho_1+\rho_2}+\frac{\rho_1|\bm u_1|^2+\rho_2|\bm u_2|^2}{3(\rho_1+\rho_2)}-\frac{|\bm u|^2}{3}\approx2.29.
\end{equation*}
}
Since one of the Gaussian distributions has a temperature $\theta_2\approx5.86 > 2 \theta$, the non-adaptive method, which uses $\bzeta=\bm u$ and $\beta=1/\sqrt{\theta}$, fails to converge. The $L^2$ error between the non-adaptive method and the initial value \eqref{eq:5BiGa} is listed in Tab.~\ref{tab:04-1}, illustrating the divergence of the non-adaptive method. 

In this simulation, the effect of the scaling adaptive method is first examined. We still set $\bzeta = \bm u$, and adjust the scaling factor $\beta$ by the adaptive scaling scheme. The parameters of the scaling adaptive method are the same as those in Tab.~\ref{tab:01-1}, with a time step length $\Delta t = 0.01$. Similarly, the initial scaling factor $\beta_0$ in \eqref{eq:beta_0} is selected by the scaling adaptive scheme. Tab.~\ref{tab:04-1} lists the initial $L^2$ errors for $N = N_0 = 10$, $15$, and $20$, using the numerical solution with $N = 24$, $N_0 = 20$ as the reference. Convergence of the scaling adaptive method is observed.

\begin{table}[!b]
    \caption{(Mixed-Gaussian problem in Sec.~\ref{sec:Eg04}) $L^2$ error between the scaling adaptive and the reference solution for $N = 10$, $15$, and $20$ at $t = 0$ and $t = 3$.}
    \label{tab:04-1}
    \centering
    \def\arraystretch{1.3}
    {\footnotesize
    \begin{tabular}{c||ccc|ccc}
     & \multicolumn{3}{c|}{$t=0$} & \multicolumn{3}{c}{$t=3$} \\
     & $N=10$ & $N=15$ & $N=20$ & $N=10$ & $N=15$ & $N=20$ \\ \hline
    scale & 3.65E$-$2 & 1.85E$-$2 & 5.52E$-$3 & 1.41E$-$6 & 2.90E$-$7 & 7.68E$-$9 \\
    non-adap & 1.83E$+$0 & 1.50E$+$1 & 1.08E$+$2 & / & / & /
    \end{tabular}
    }
\end{table}

\begin{figure}
    \centering
    \subfigure[$L^2$ error\label{fig:04-1a}]{
        \includegraphics[height=7.4\baselineskip]{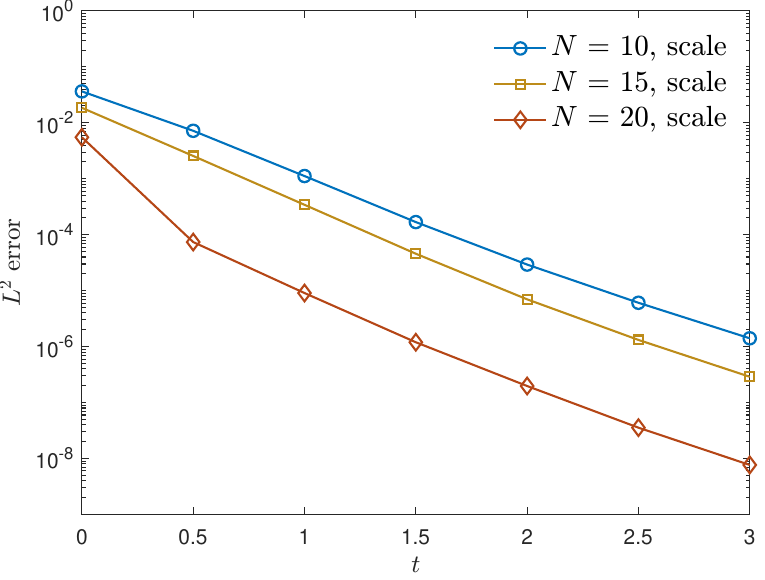}
    }\hspace{-0.7em}
    \subfigure[scaling factor $\beta$\label{fig:04-1b}]{
        \includegraphics[height=7.4\baselineskip]{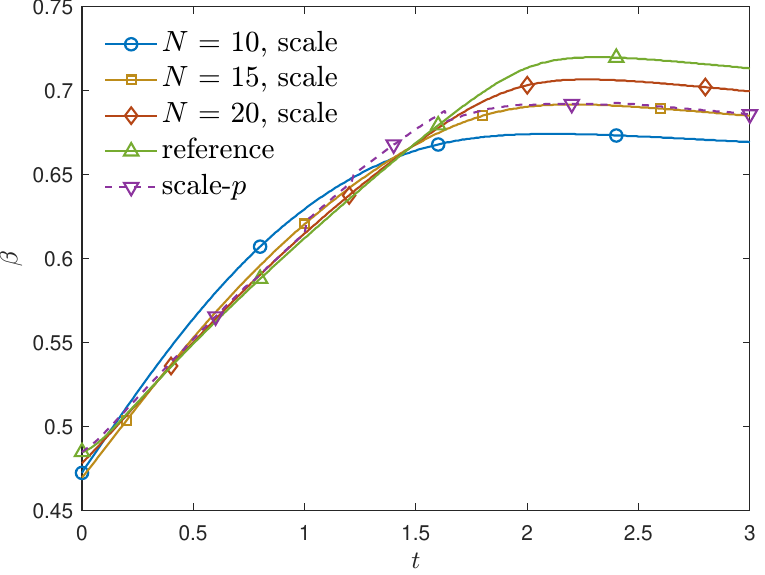}
    }\hspace{-0.7em}
    \subfigure[$L^2$ error and $N$\label{fig:04-1c}]{
        \includegraphics[height=7.4\baselineskip]{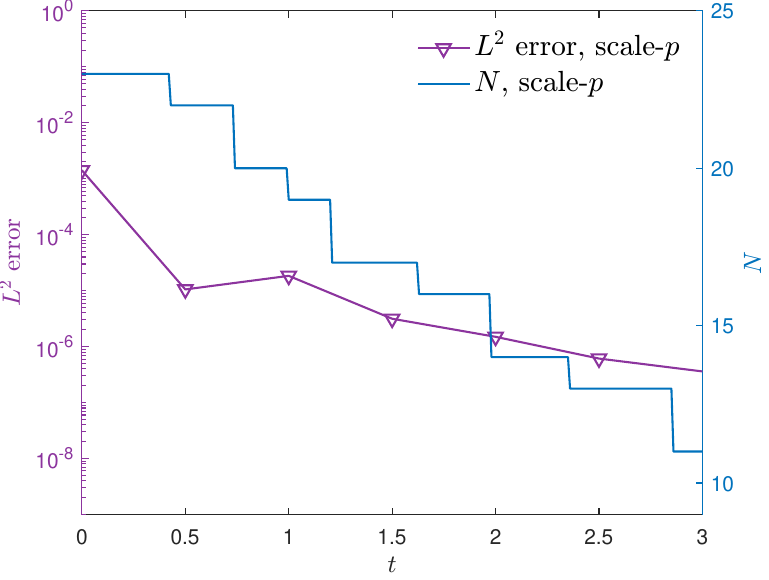}
    }
    
    \caption{(Mixed-Gaussian problem in Sec.~\ref{sec:Eg04}) Performance of the scaling adaptive and scale-$p$-adaptive methods with different expansion order $N$. (a) Evolution of the $L^2$ error between the scaling adaptive and reference solution. (b) Evolution of the scaling factor $\beta$. (c) Evolution of the $L^2$ error (purple) and expansion order $N$ (blue) for the scale-$p$-adaptive method.}
    \label{fig:04-1}
\end{figure}

\begin{figure}
    \centering
    \subfigure[$t=0$\label{fig:04-2a}]{
        \includegraphics[height=7.4\baselineskip]{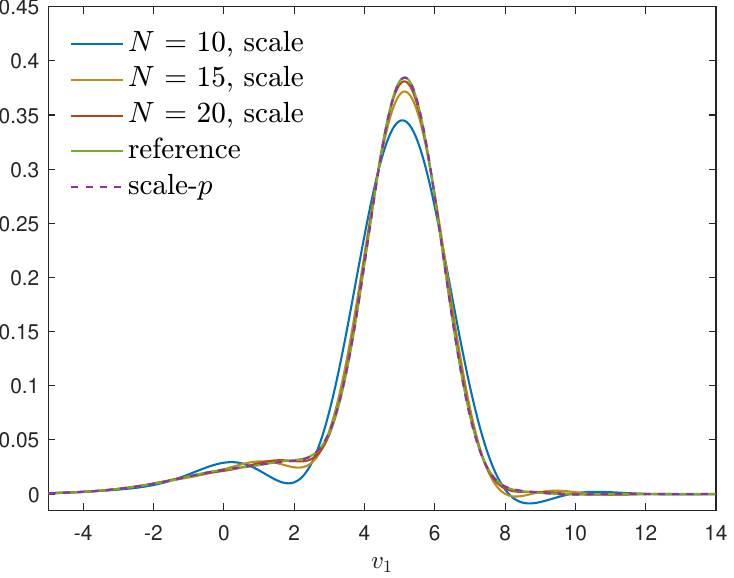}
    }\hspace{-0.7em}
    \subfigure[$t=0.25$\label{fig:04-2b}]{
        \includegraphics[height=7.4\baselineskip]{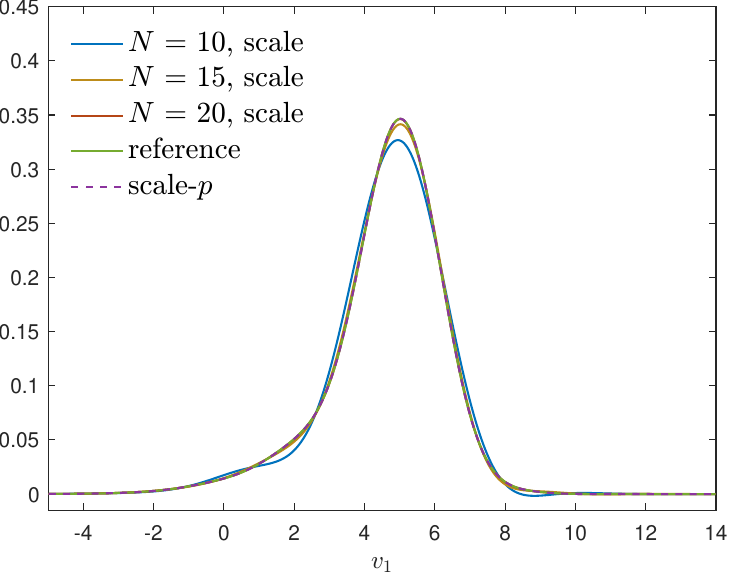}
    }\hspace{-0.7em}
    \subfigure[$t=3$\label{fig:04-2c}]{
        \includegraphics[height=7.4\baselineskip]{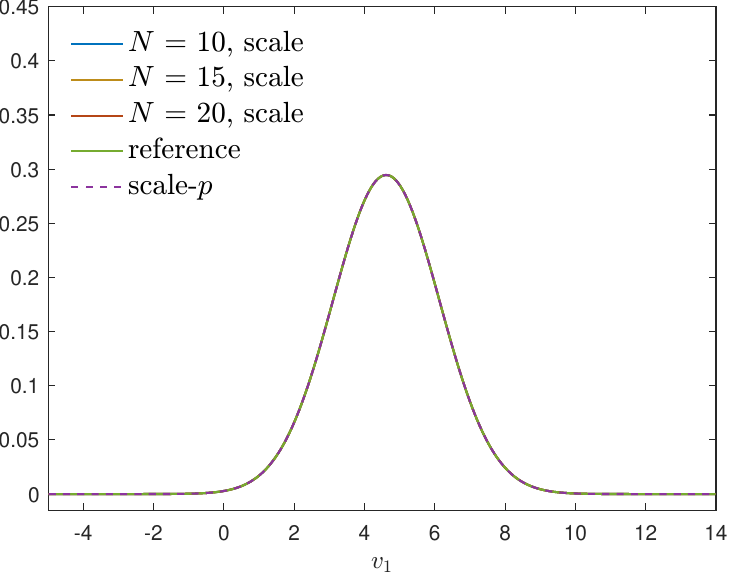}
    }
    
    \caption{(Mixed-Gaussian problem in Sec.~\ref{sec:Eg04}) Marginal distribution functions at different times. Results from the scaling adaptive and scale-$p$-adaptive methods are compared with the reference solution. (a) $t = 0$. (b) $t = 0.25$. (c) $t = 3$.}
    \label{fig:04-2}
\end{figure}

The evolution of the $L^2$ error and the scaling factor $\beta$ are shown in Fig.~\ref{fig:04-1a} and \ref{fig:04-1b}, respectively, and the $L^2$ errors at $t = 3$ are also given in Tab.~\ref{tab:04-1}. Results indicate that the scaling adaptive method converges as $N$ increases, and that $\beta$ is gradually adjusted to $1/\sqrt{\theta}\approx 0.66$, which is consistent with the expectation as the distribution function $f$ evolves toward the Maxwellian. Fig.~\ref{fig:04-2} compares the numerical solutions for $N = 10$, $15$, and $20$ with the reference solution at $t = 0$, $0.25$, and $3$. For $N = 10$ and $15$, noticeable discrepancies are observed at $t = 0$ and $0.25$, whereas for $N = 20$, the numerical solution is nearly identical to the reference. By $t = 3$, all these numerical solutions and the reference solution are on top of each other.

The rapid decrease in the $L^2$ error indicates that the expansion order $N$ can be reduced during the later stages of the computation. To achieve this and lower the computational cost, the scale-$p$-adaptive method is tested here. The initial expansion order is set to $N = 23$, and the parameters for the $p$-adaptive algorithm are listed in Tab.~\ref{tab:04-2}. Fig.~\ref{fig:04-1c} shows the evolution of both the $L^2$ error and $N$. As the solution approaches equilibrium, the expansion order decreases gradually from $N = 23$ to $N = 11$, while the $L^2$ error continues to decline, remaining within the order of $10^{-6}$ at $t = 3$. This demonstrates the effectiveness of the scale-$p$-adaptive method. The scaling factor for this method, shown in Fig.~\ref{fig:04-1b}, also approaches $1/\sqrt{\theta}$ as expected.

The computational cost, as summarized in Tab.~\ref{tab:04-3}, exhibits behavior consistent with the previous two examples. Furthermore, the scale-$p$-adaptive method requires only half the computational time of the scaling adaptive method with $N = 20$, confirming its high efficiency.

\begin{table}[hptb]
    \caption{(Mixed-Gaussian problem in Sec.~\ref{sec:Eg04}) Parameters utilized in the $p$-adaptive algorithm.}
    \label{tab:04-2}
    \centering
    \def\arraystretch{1.3}
    {\footnotesize
    \begin{tabular}{c|cccccccc}
    parameter & $N_{\min}$ & $N_{\max}$ & $\Delta N$ & $\eta_{l}^{(p)}$ & $\eta_{h}^{(p)}$ & $\mathcal{F}_{l,0}$ & $\mathcal{F}_{h,0}$ & $\eta_0^{(p)}$ \\ \hline
    usage & Alg.~\ref{alg:ada-p} & Alg.~\ref{alg:ada-p} & Alg.~\ref{alg:ada-p} & \eqref{eq:thre} & \eqref{eq:thre} & \eqref{eq:thre} & \eqref{eq:thre} & \eqref{eq:update_F_new} \\ 
    value & 2 & 24 & 1 & 0.3 & 1.5 & $10^{-9}$ & $10^{-6}$ & 0.01
    \end{tabular}
    }
\end{table}

\begin{table}[hptb]
    \caption{(Mixed-Gaussian problem in Sec.~\ref{sec:Eg04}) Average CPU time per time step (in seconds). Here, $T_{\rm non}$ and $T_{\rm adap}$ refer to the CPU time of the non-adaptive and scaling adaptive methods, respectively. $T_{\rm ind}$ refers to the CPU time spent on the adaptive adjustment step, which is a component of $T_{\rm adap}$.}
    \label{tab:04-3}
    \centering
    \def\arraystretch{1.3}
    {\footnotesize
    \begin{tabular}{cc||ccc}
     &  & $T_{\rm adap}$ & $T_{\rm ind}$ & $T_{\rm ind}/T_{\rm adap}$ \\ \hline
    \multirow{3}{*}{scale} & $N=10$ & 4.64E$-$2 & 4.88E$-$5 & 0.105\% \\
     & $N=15$ & 1.05E$+$0 & 1.72E$-$4 & 0.016\% \\
     & $N=20$ & 9.48E$+$0 & 4.56E$-$4 & 0.005\% \\ \hline
    \multicolumn{2}{c||}{scale-$p$} & 4.40E$+$0 & 4.35E$-$4 & 0.009\%
    \end{tabular}
    }
\end{table}

\subsection{Spatially non-homogeneous case}
In this subsection, the adaptive methods are tested on four non-homogeneous problems. The IPL model with $\eta = 10$ is employed, with the order $N_0 = 15$ in \eqref{eq:5Qhat}. The time step length is determined by the CFL number as
{\small
\begin{equation*}
    {\rm CFL}=\sum_{d=1}^{D}\frac{v_{\max,d}}{\Delta x_d/\Delta t},
\end{equation*}
}
where $\Delta x_d$ denotes the mesh size, and $v_{\max, d}$ is the maximum characteristic velocity in the $d$-th dimension ($d = 1, \dots, D$). We refer \cite{Hu2020} for the detailed discussion of $v_{\max, d}$. Since all the initial values in this section are symmetric in velocity space, $\bzeta$ is set as $(0,0,0)$. A periodic boundary condition is applied to all problems.

\subsubsection{1D: Maxwellian with perturbation on density}
\label{sec:Eg11}
First, the performance of the scaling adaptive scheme is tested on a spatially one-dimensional wave problem, where the initial condition is set as
{\small
\begin{gather}
    f(0,x,\bm v)=\frac{\rho(x)}{(2\pi\theta(x))^{3/2}}\exp\left(-\frac{|\bm{v}|^2}{2\theta(x)}\right),\label{eq:13Maxw}\\
    \rho(x)=1+0.3\sin(\pi x),\qquad \theta(x)=1, 
    \qquad x\in[0,2].
\end{gather}
}
Similar problems have been tested in the literature \cite{Xiao2020}. The Knudsen number is set as $\varepsilon = 5$, the expansion order of the Fourier method in the spatial domain is $M = 16$, and ${\rm CFL} = 0.5$. The parameters utilized in the scaling adaptive scheme are the same as Tab.~\ref{tab:01-1}. The reference solution is obtained by the non-adaptive method with fixed $\beta \equiv 1$ and $N = 50$.

\begin{figure}
    \centering
    \subfigure[$\rho$\label{fig:11-1a}]{
        \includegraphics[height=7.4\baselineskip]{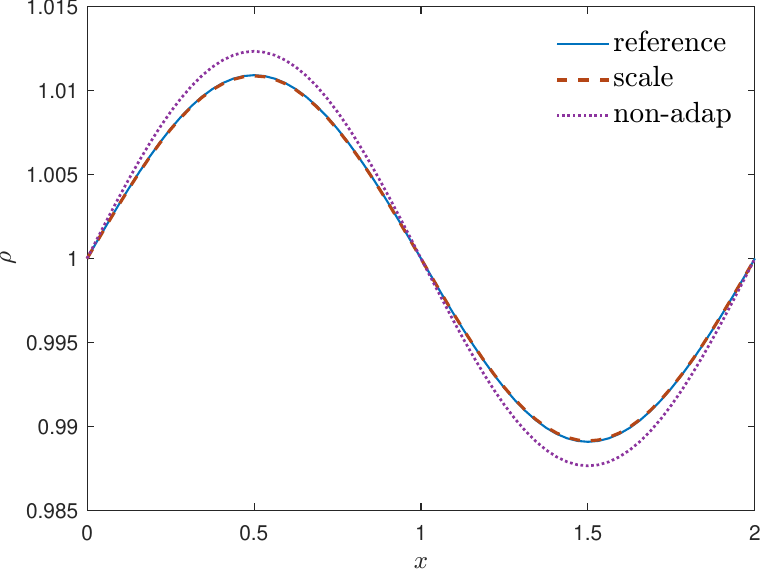}
    }\hspace{-0.7em}
    \subfigure[$\theta$\label{fig:11-1b}]{
        \includegraphics[height=7.4\baselineskip]{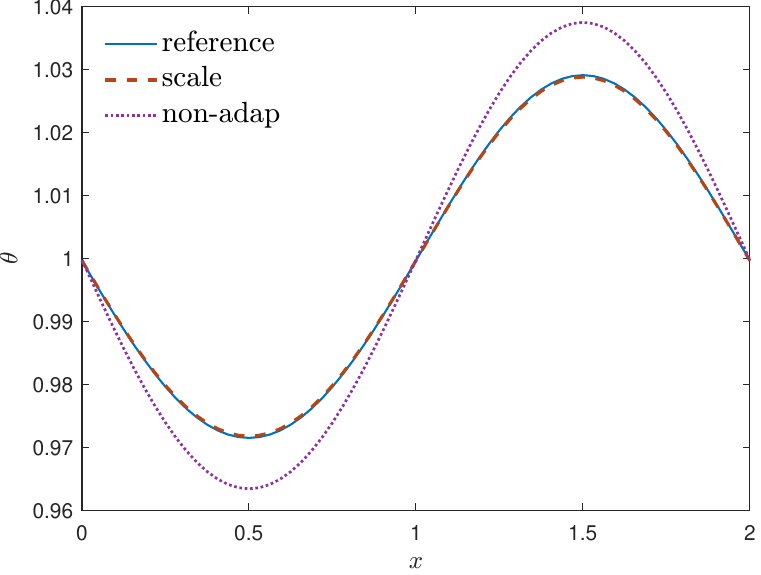}
    }\hspace{-0.7em}
    \subfigure[$L^2$ error of $\rho$ and $\theta$\label{fig:11-1c}]{
        \includegraphics[height=7.4\baselineskip]{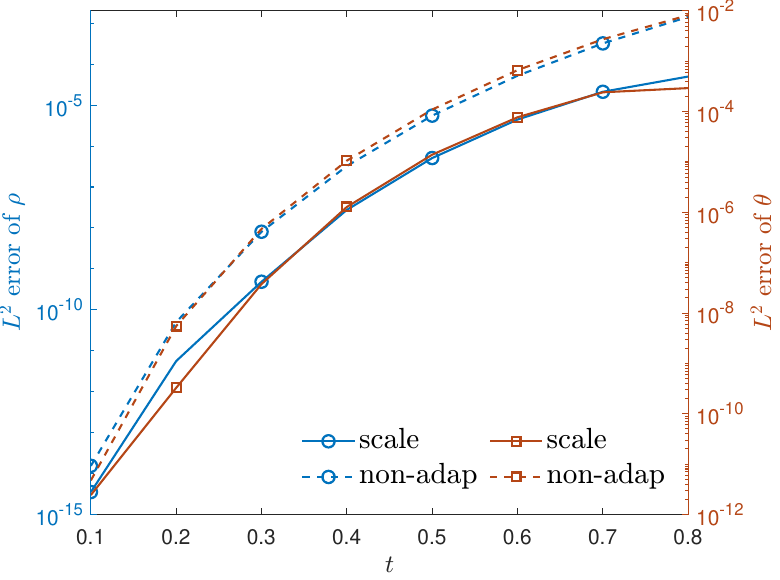}
    }
    
    \caption{(1D: Maxwellian perturbation problem in Sec.~\ref{sec:Eg11}) Macroscopic variables of the scaling adaptive and non-adaptive methods at $t=0.8$ with $N = 6$. (a) Density $\rho$. (b) Temperature $\theta$. (c) Evolution of the $L^2$ error between the numerical and the reference macroscopic density (blue) and temperature (red).}
    \label{fig:11-1}
\end{figure}

For this case, the initial scaling factor is set to $\beta_0 = 1$ in the scaling adaptive method, since the initial condition \eqref{eq:13Maxw} is a standard Maxwellian distribution function with the same velocity and temperature. The density and temperature at $t = 0.8$ obtained by the scaling adaptive and non-adaptive method with $N = 6$ are plotted in Fig.~\ref{fig:11-1a} and \ref{fig:11-1b}, respectively. The scaling adaptive method produces results that closely match the reference solution, whereas the non-adaptive method with the same expansion order shows significant deviation. The evolution of the $L^2$ error for the density and temperature obtained by the scaling adaptive and the non-adaptive method is shown in Fig.~\ref{fig:11-1c}. It reveals that the $L^2$ error obtained by the scaling adaptive method is much smaller than that of the non-adaptive method, indicating the high accuracy of the scaling adaptive method. 

The CPU times of the scaling adaptive and non-adaptive methods are listed in Tab.~\ref{tab:11-1}. The behavior of the scaling adaptive method is similar to the non-adaptive method, as the computational time for the two methods is almost the same, and the ratio of calculating the indicator is quite small. 

\begin{table}[hptb]
    \caption{(Non-homogeneous problems) The average computational time (wall time in seconds) for each time step. Here, $T_{\rm non}$ and $T_{\rm adap}$ refer to the computational time of the non-adaptive and the adaptive method. $T_{\rm ind}$ refers to the computational time of the adjustment step in $T_{\rm adap}$. Tests in Sec.~\ref{sec:Eg21} and \ref{sec:Eg31} are run with 32 threads.
    }
    \label{tab:11-1}
    \centering
    \def\arraystretch{1.3}
    {\footnotesize
    \begin{tabular}{c||cccc}
     & $T_{\rm non}$ & $T_{\rm adap}$ & $T_{\rm ind}$ & $T_{\rm ind}/T_{\rm adap}$ \\ \hline
    Sec.~\ref{sec:Eg11} & 1.09E$-$2 & 1.12E$-$2 & 2.56E$-$4 & 2.286\% \\
    Sec.~\ref{sec:Eg14} & 1.17E$+$1 & 0.77E$+$1 & 2.04E$-$5 & 0.000\% \\
    Sec.~\ref{sec:Eg21} & 2.40E$+$0 & 6.22E$-$1 & 1.85E$-$3 & 0.297\% \\
    Sec.~\ref{sec:Eg31} & 1.92E$+$0 & 1.92E$+$0 & 4.56E$-$3 & 0.238\%
    \end{tabular}
    }
\end{table}

\subsubsection{1D: Time-varying regime problem}\label{sec:Eg14}
\begin{figure}
    \centering
    \subfigure[Knudsen number $\varepsilon$ and $N$\label{fig:14-1a}]{
        \includegraphics[height=7.33\baselineskip]{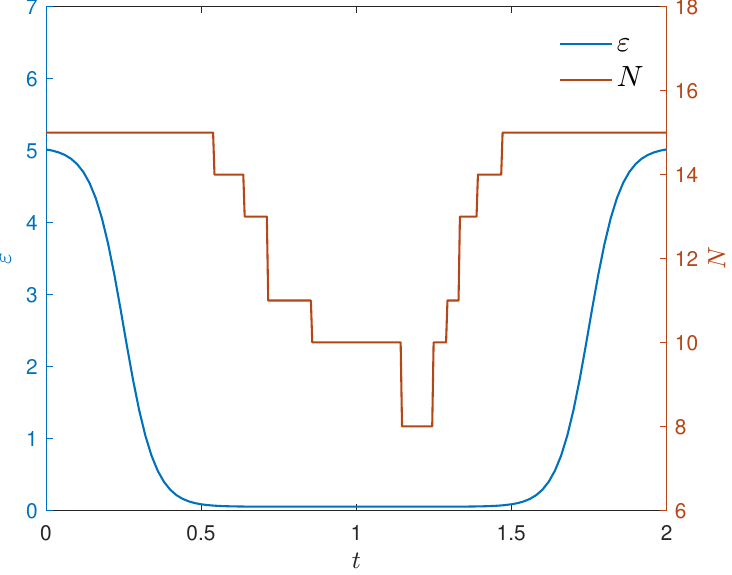}
    }\hspace{-0.7em}
    \subfigure[frequency indicator $\mF$\label{fig:14-1b}]{
        \includegraphics[height=7.4\baselineskip]{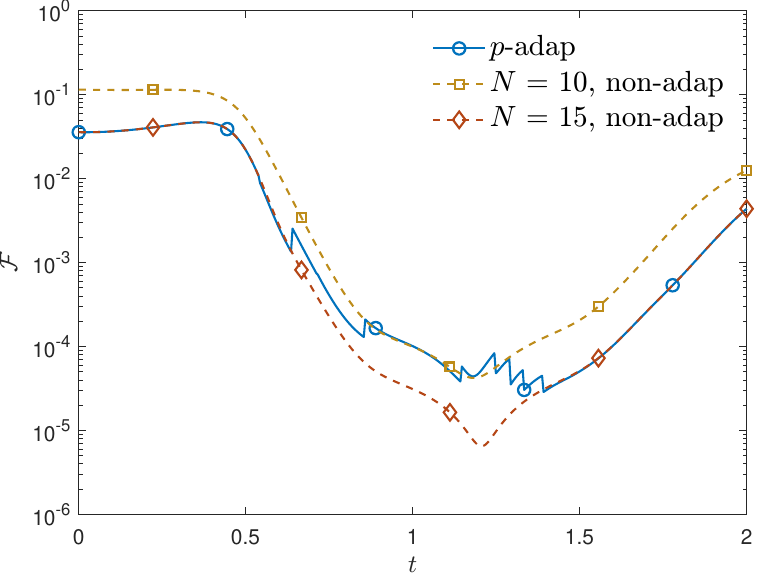}
    }\hspace{-0.7em}
    \subfigure[$L^2$ error of $f$\label{fig:14-1c}]{
        \includegraphics[height=7.4\baselineskip]{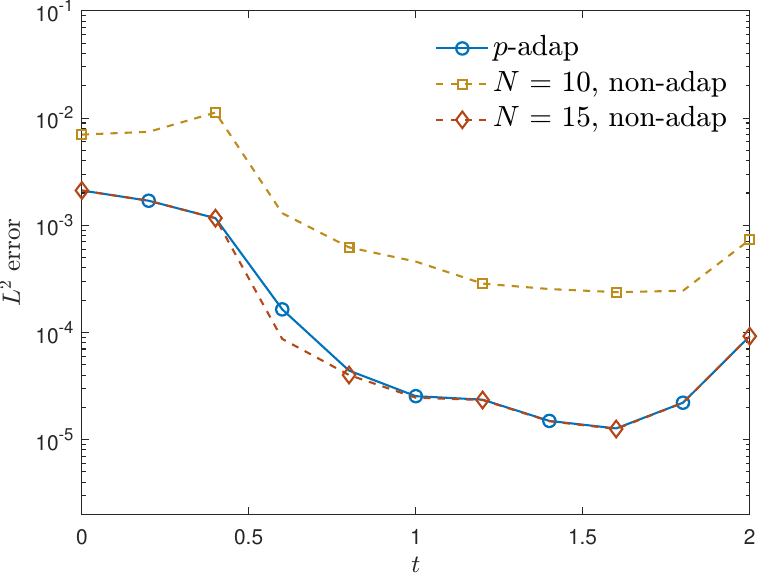}
    }

    \caption{(1D: Time-varying regime problem in Sec.~\ref{sec:Eg14}) Comparison between the $p$-adaptive and non-adaptive methods. (a) Knudsen number $\varepsilon$ setting (blue) and the evolution of the expansion order $N$ (red). (b) Evolution of the frequency indicator $\mF$. (c) Evolution of the $L^2$ error between the numerical and the reference solution $f$.}
    \label{fig:14-1}
\end{figure}

The $p$-adaptive method is tested in this section on a time-varying regime problem, where the Knudsen number is time-dependent as in 
{\small
\begin{equation*}
    \varepsilon(t)=5.05+2.5(\tanh(10(0.25-t))+\tanh(10(t-1.75))).
\end{equation*}
}
As illustrated in Fig.~\ref{fig:14-1a}, the initial Knudsen number is $\mathcal{O}(1)$, then reduced to $\mathcal{O}(0.01)$, and finally back to $\mathcal{O}(1)$. In the intermediate interval $t\in[0.5,1.5]$, the solution $f$ will approach the local equilibrium, which may require a smaller $N$. The initial condition of this problem is a combination of two Maxwellian distributions as
{\small
\begin{align*}
    &f(0,x,\bm v)=\frac{\rho(x)}{(2\pi\theta(x))^{3/2}}\exp\left(-\frac{|\bm{v}-\bm{u}|^2}{2\theta(x)}\right)+\frac{\rho(x)}{(2\pi\theta(x))^{3/2}}\exp\left(-\frac{|\bm{v}+\bm{u}|^2}{2\theta(x)}\right),\\ 
    &\rho(x)=\frac{1+0.1\sin(\pi x)}{2},\quad \bm{u}=(\sqrt{3/2},0,0),\quad \theta(x)=\frac{1-0.1\sin(\pi x)}{2},\quad x\in[0,2].
\end{align*}
}
In the simulation, the expansion order of the Fourier method in the spatial space is $M = 16$, and the $\rm CFL$ number is set as $\rm CFL = 0.5$. In the $p$-adaptive method, $N_{\max} = 15$, and other parameters are the same as Tab.~\ref{tab:04-2}. The $p$-adaptive method with initial expansion order $N=15$ is tested. The reference solution is obtained by the non-adaptive method with fixed $\beta \equiv 1$ and $N = 50$.

\begin{figure}
    \centering
    \subfigure[$\rho$\label{fig:14-2a}]{
        \includegraphics[height=7.4\baselineskip]{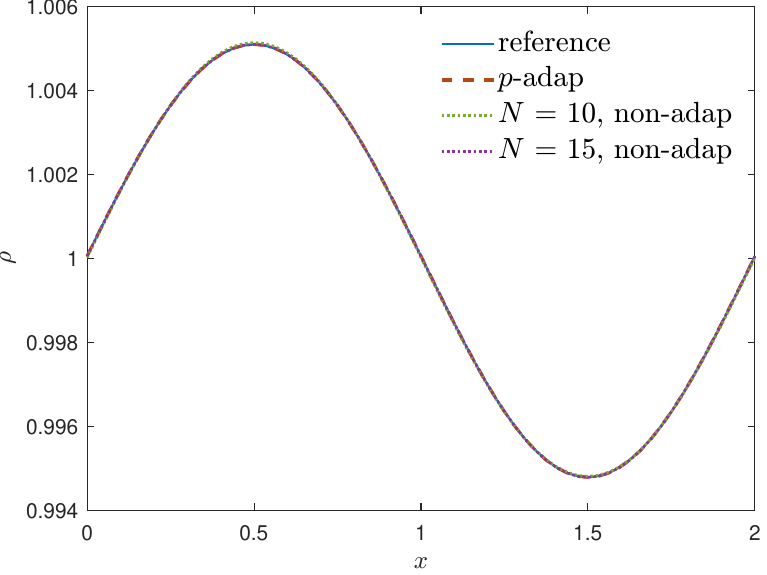}
    }\hspace{-0.7em}
    \subfigure[$\theta$\label{fig:14-2b}]{
        \includegraphics[height=7.4\baselineskip]{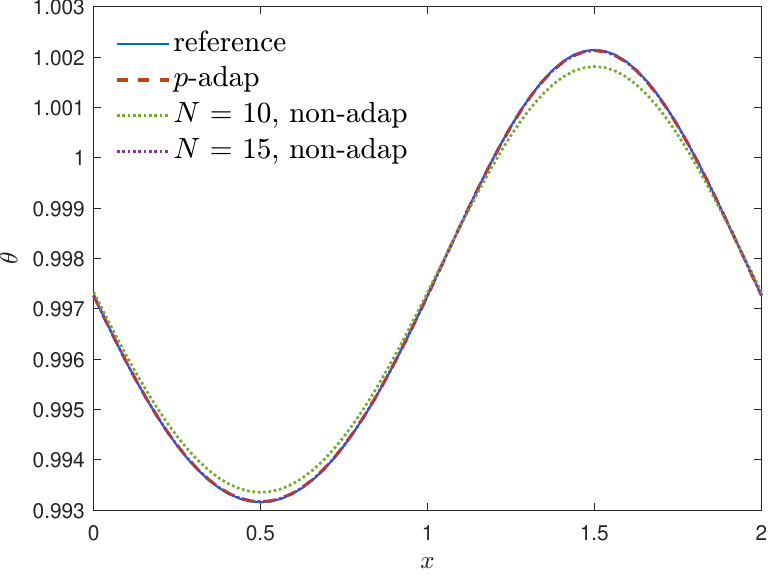}
    }\hspace{-0.7em}
    \subfigure[$q_1$\label{fig:14-2c}]{
        \includegraphics[height=7.7\baselineskip]{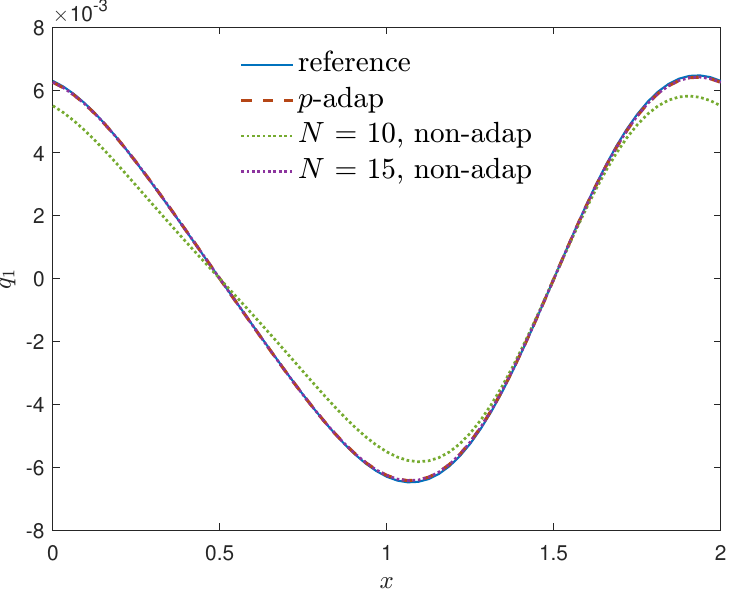}
    }

    \caption{(1D: Time-varying regime problem in Sec.~\ref{sec:Eg14}) Macroscopic variables of the $p$-adaptive and non-adaptive methods at $t=2$ with different expansion order $N$. (a) Density $\rho$. (b) Temperature $\theta$. (c) Heat flux $q_1$.}
    \label{fig:14-2}
\end{figure}

The numerical solution of the density $\rho$, the temperature $\theta$, and the heat flux $q_1$ at $t = 2$ is plotted in Fig.~\ref{fig:14-2}, where the reference solution and the numerical solution with the non-adaptive method with $N = 10$ and $N = 15$ are also plotted. It shows that for the density $\rho$, all these solution matches well with the reference solution. However, for the temperature $\theta$ and heat flux $q_1$, the numerical solution obtained by the $p$-adaptive method and that by the non-adaptive method with $N = 15$ are nearly the same as the reference solution, while that by the non-adaptive method with $N = 10$ has an evident difference from the reference solution. The evolution of the indicator $\mF$ and the $L^2$ error for the distribution function $f$ is shown in Fig.~\ref{fig:14-1b} and \ref{fig:14-1c}, respectively. The error of the $p$-adaptive and the non-adaptive method with $N = 15$ is almost the same, and is much smaller than that of the non-adaptive method with $N = 10$. This means that the solution of this $p$-adaptive method is similar to the non-adaptive method with $N = 15$, and much better than that with $N = 10$. The evolution of the expansion number $N$ for the $p$-adaptive method is illustrated in Fig.~\ref{fig:14-1a}. The tendency of $N$ is similar to the Knudsen number, and the minimum decrease is $N = 8$ when $\varepsilon$ is small. 

Therefore, a significant amount of computation is saved. The CPU time for the $p$-adaptive method and the non-adaptive method with expansion order $N=15$ are listed in Tab.~\ref{tab:11-1}. With almost identical error, the $p$-adaptive method saved one-third of the computational time. The proportion of adaptive adjustment time in the total method is less than $10^{-5}$, which confirms the efficiency of the $p$-adaptive algorithm.

\subsubsection{2D: Taylor-Green vortex}
\label{sec:Eg21}
The spatially 2-dimensional Taylor-Green vortex problem is studied. Its initial condition is a local Maxwellian with the density, macroscopic velocity, and temperature as below 
{\small
\begin{equation*}
\begin{gathered}
    \rho(x,y)=1, \quad \theta(x,y)=\theta_0-{u_{0}^2}(\cos(4\pi x)+\cos(4\pi y))/4, \quad \theta_0 = 1,\\
    \bm{u}(x,y)=(-u_0\cos(2\pi x)\sin(2\pi y),u_0\sin(2\pi x)\cos(2\pi y),0), \quad u_0 = 1.2.
\end{gathered}
\end{equation*}
}
with the spatial domain $\Omega=[0,1]^2$. This Taylor-Green vortex problem has an analytic solution for the incompressible Navier-Stokes equation \cite{Chorin1968}, where the macroscopic velocity decays to $0$ as time increases. For the Boltzmann equation simulated here, the Knudsen number is set as $\varepsilon=0.05$, and the distribution function is expected to converge to the uniform Maxwellian with the final density, macroscopic velocity, and temperature as below
{\small
\begin{equation}\label{eq:Eg31end}
    \rho(x,y)=1,\qquad \bm{u}(x,y)=(0,0,0),\qquad \theta(x,y)=\theta_0+u_0^2/6=1.24.
\end{equation}
}

In the simulation, the expansion order of the Fourier method in spatial space is set as $M_x = M_y = 32$ with ${\rm CFL}=0.5$. The parameters utilized in the scale-$p$-adaptive method are the same as Tab.~\ref{tab:01-1} and \ref{tab:04-2}, except for $\eta_{l}^{(p)}=0.1$. The initial expansion order for the scale-$p$-adaptive method is set as $N = 10$. The reference solution is obtained by the non-adaptive method with fixed $\beta\equiv1/\sqrt{\theta_0}=1$ and $N=50$.

\begin{figure}
    \centering
    \subfigure[$u_1$\label{fig:21-2a}]{
        \includegraphics[width=0.32\linewidth]{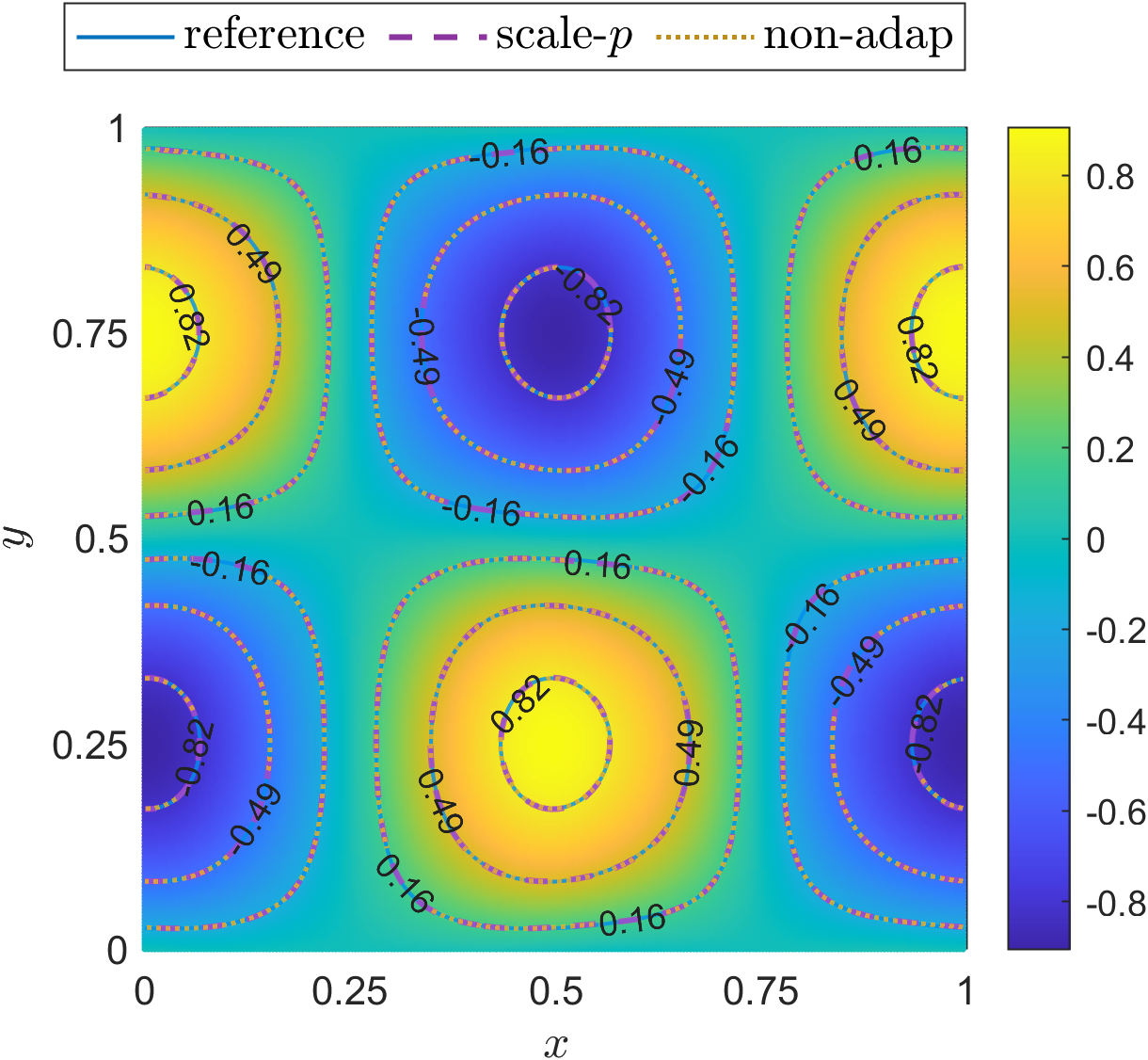}
    }\hspace{-0.7em}
    \subfigure[$u_2$\label{fig:21-2b}]{
        \includegraphics[width=0.32\linewidth]{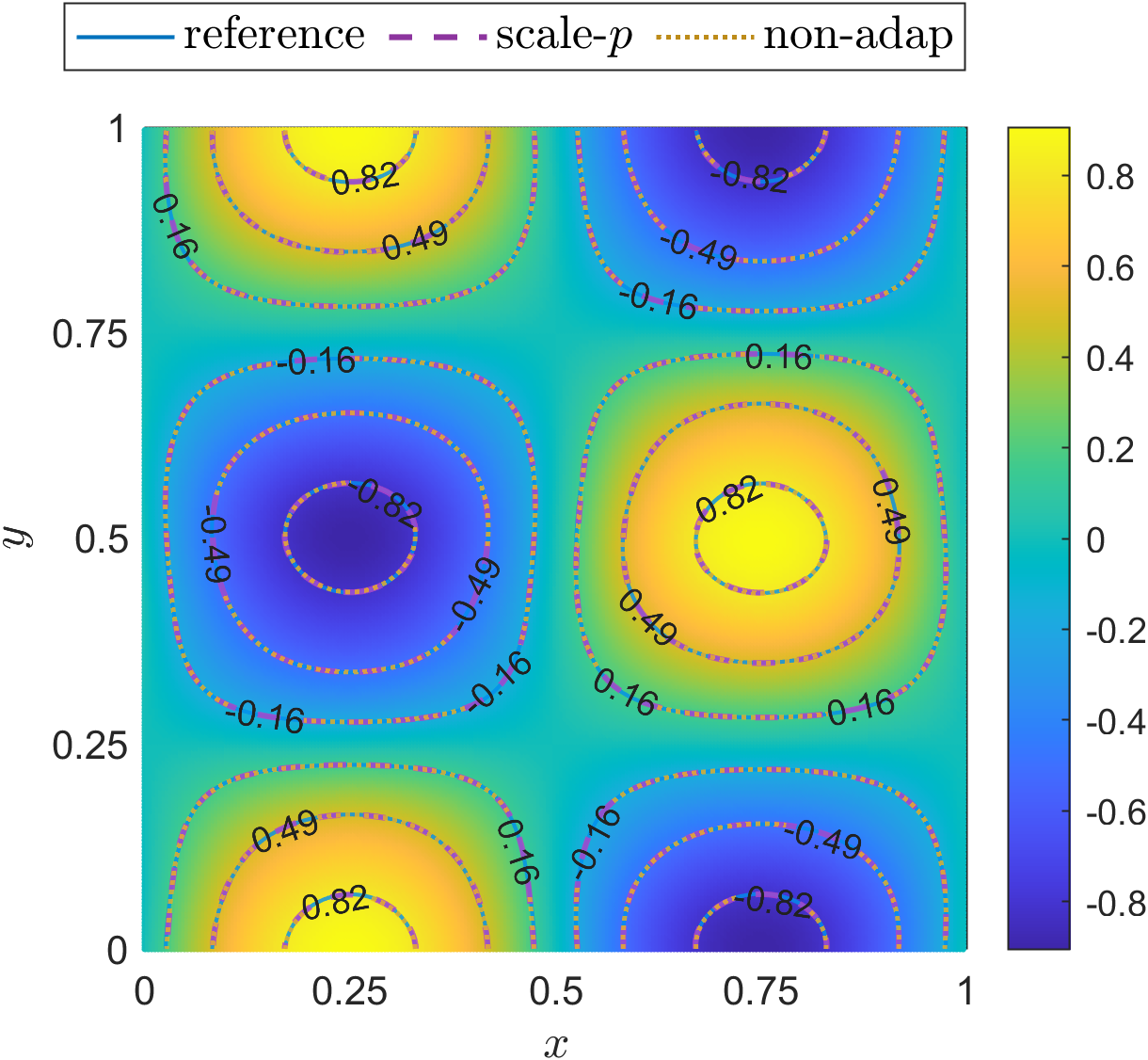}
    }\hspace{-0.7em}
    \subfigure[$\theta$\label{fig:21-2c}]{
        \includegraphics[width=0.32\linewidth]{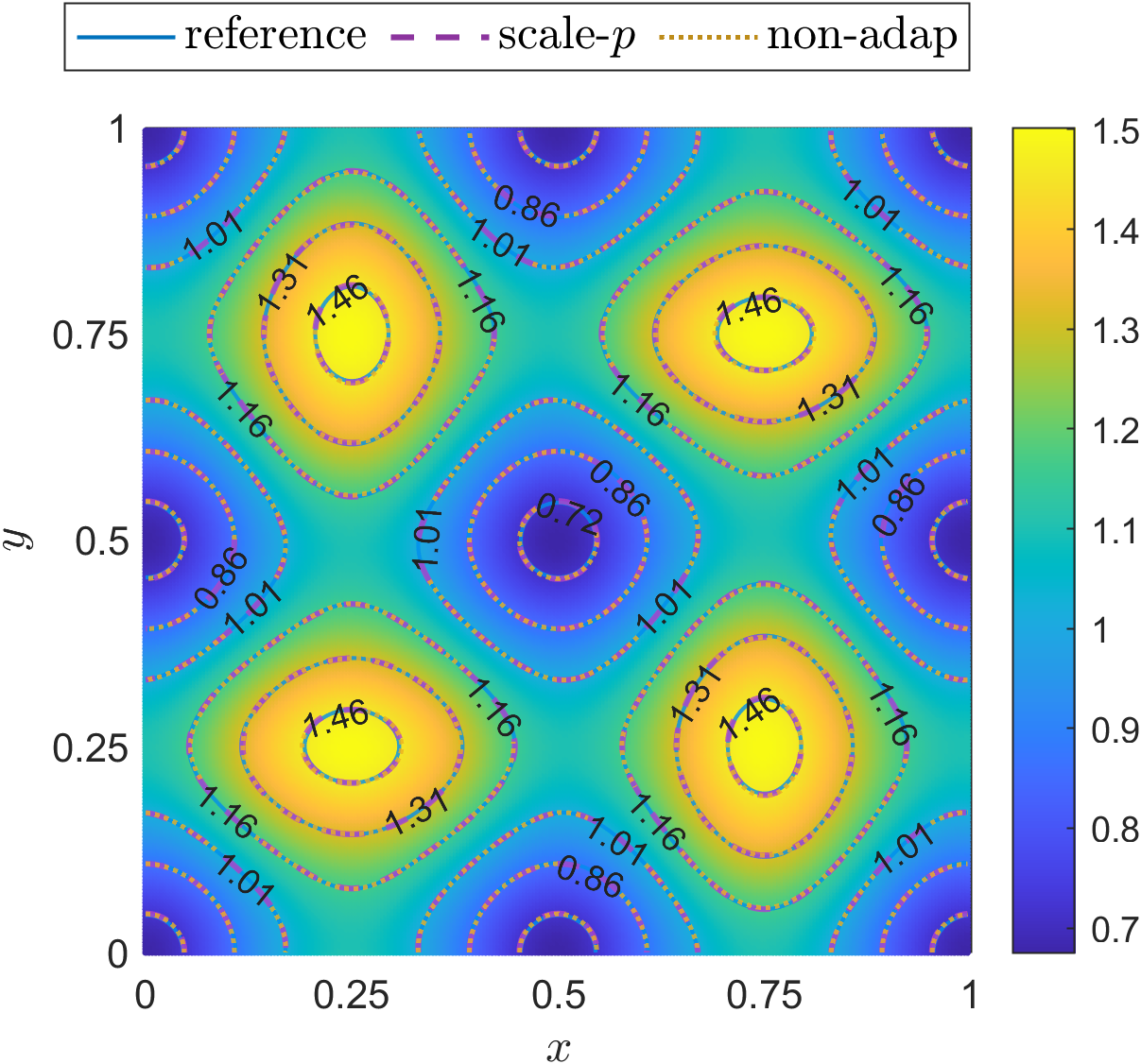}
    }
    
    \caption{(2D: Taylor-Green vortex in Sec.~\ref{sec:Eg21}) Macroscopic variables of the scale-$p$-adaptive and non-adaptive method with $N=10$, $t=0.1$. (a) Velocity $u_1$. (b) Velocity $u_2$. (c) Temperature.}
    \label{fig:21-2}
\end{figure}

\begin{figure}
    \centering
    \subfigure[$\theta(x,0.25)$\label{fig:21-1a}]{
        \includegraphics[height=7.5\baselineskip]{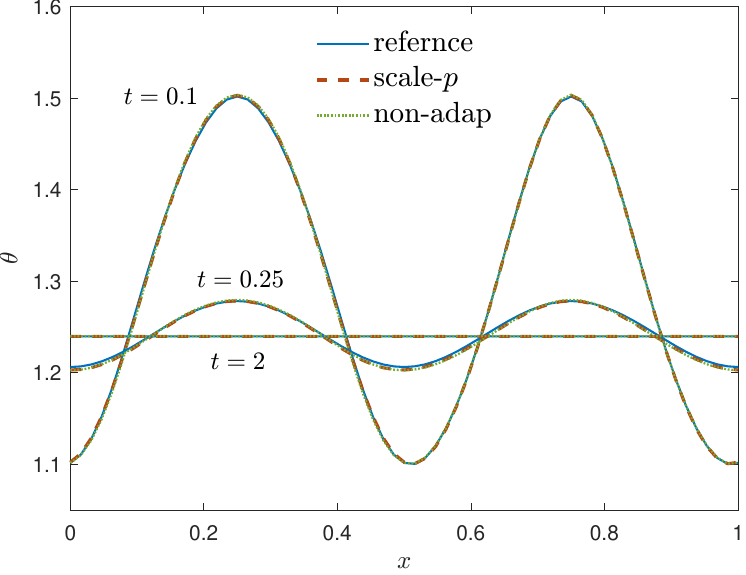}
    }\hspace{-0.7em}
    \subfigure[$L^2$ error and $N$\label{fig:21-1b}]{
        \includegraphics[height=7.5\baselineskip]{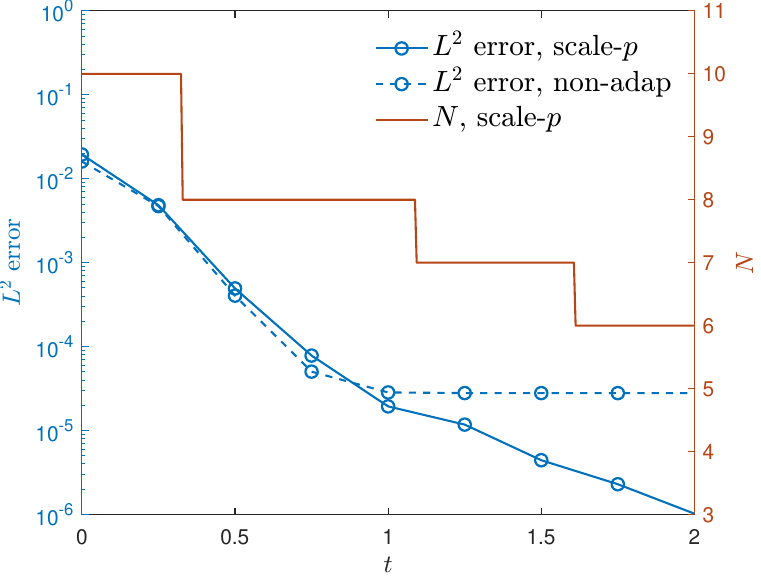}
    }\hspace{-0.7em}
    \subfigure[$\mF$ and $\beta$\label{fig:21-1c}]{
        \includegraphics[height=7.5\baselineskip]{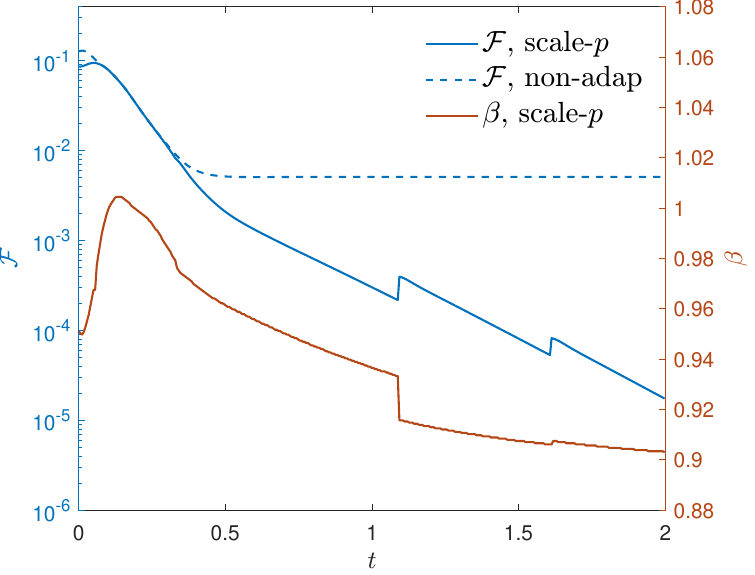}
    }
    
    \caption{(2D: Taylor-Green vortex in Sec.~\ref{sec:Eg21}) Comparison between the scale-$p$-adaptive and non-adaptive methods with $N=10$. (a) Temperature $\theta$ along the line $y=0.25$ at times $t=0.1,0.25,2$. (b) Evolution of the $L^2$ error between the numerical and the reference solution (blue) and the expansion order $N$ (red) of the scale-$p$-adaptive method. (c) Evolution of the frequency indicator $\mF$ (blue) and the scaling factor $\beta$ (red) of the scale-$p$-adaptive method.}
    \label{fig:21-1}
\end{figure}

The numerical solution of the macroscopic velocity $u_1$ and $u_2$, the temperature $\theta$ at $t = 0.1$ is plotted in Fig.~\ref{fig:21-2}, where the numerical solution by the non-adaptive Hermite method with $N \equiv 10$ and the reference solution are plotted. It shows that for all the macroscopic variables, the three solutions all match well with each other. Moreover, the temperature along the line $y = 0.25$ at $t = 0.1, 0.25$, and $2$ is plotted in Fig.~\ref{fig:21-1a}, where the consistency of the numerical solution by the scale-$p$-adaptive method and the reference is verified more clearly. 

The evolution of $L^2$ error of the distribution function and the variation of the expansion order for the adaptive method, as well as the $L^2$ error for the non-adaptive method, are shown in Fig.~\ref{fig:21-1b}. It shows that the expansion order is reducing gradually, and the minimum order is $N = 6$. Even though the $L^2$ error is nearly the same as the non-adaptive method at the beginning, it becomes smaller when $t > 1$. The evolution of the frequency indicator $\mF$ and scaling factor $\beta$ is illustrated in Fig.~\ref{fig:21-1c}. The scaling factor is near $1$ at the early stage and gradually adjusted to $1/\sqrt{\theta_0+u_0^2/6}\approx 0.8980$, which matches the theoretical scaling factor of the final equilibrium with macroscopic variables \eqref{eq:Eg31end}.


The computational time for the scale-$p$-adaptive method and the non-adaptive method with expansion order $N=10$ is listed in Tab.~\ref{tab:11-1}. Both tests utilize OpenMP with 32 threads, and the wall time is presented. The computational time of the scale-$p$-adaptive method is approximately only $26\%$ of that of the non-adaptive method, demonstrating the effectiveness of the adaptive algorithm.

\subsubsection{3D: Hexa-Gaussian with perturbation on density}\label{sec:Eg31}
To further illustrate the high efficiency of this adaptive Hermite method, a spatially three-dimensional problem is considered in this section. The initial condition is a Hexa-Gaussian with a perturbation on the density
{\small
\begin{gather*}
    f(0,x,y,z,\bm{v})= \sum_{i = 1}^6 \frac{\rho(x,y,z)}{(2\pi\theta)^{3/2}}\exp\left(-\frac{|\bm{v}-\bm{u}_i|^2}{2\theta}\right),\\
    \rho(x,y,z) = \frac{1+0.3\lambda(x,y,z)}{6},\qquad \theta = 1/2, \qquad (x,y,z)\in[0,2]^3,\\ 
    \bm{u}_1=-\bm{u}_2=(\sqrt{3/2},0,0),
    \quad \bm{u}_3=-\bm{u}_4=(0,\sqrt{3/2},0), \quad \bm{u}_5=-\bm{u}_6=(0,0,\sqrt{3/2}),\\
    \lambda(x,y,z)=(\sin(\pi x)-0.5\sin(2\pi x)+0.125\sin(3\pi x))(\sin(\pi y)-0.5\sin(2\pi y))\sin(\pi z).
\end{gather*}
}
Due to the high computational cost of the 3D problem, a small expansion number $N = 8$ is utilized, and only the scaling adaptive method is tested. In the simulation, the spatial space is discretized by the Fourier spectral method with $M_x=M_y=M_z=16$ and CFL number ${\rm CFL}=0.5$. The parameters utilized in the scaling adaptive method are the same as Tab.~\ref{tab:01-1}. The reference solution is obtained by the non-adaptive method with fixed $\beta\equiv1$ and $N=50$.

\begin{figure}
    \centering
    \subfigure[$\rho$, $z=0.5$\label{fig:31-1a}]{
        \includegraphics[width=0.32\linewidth]{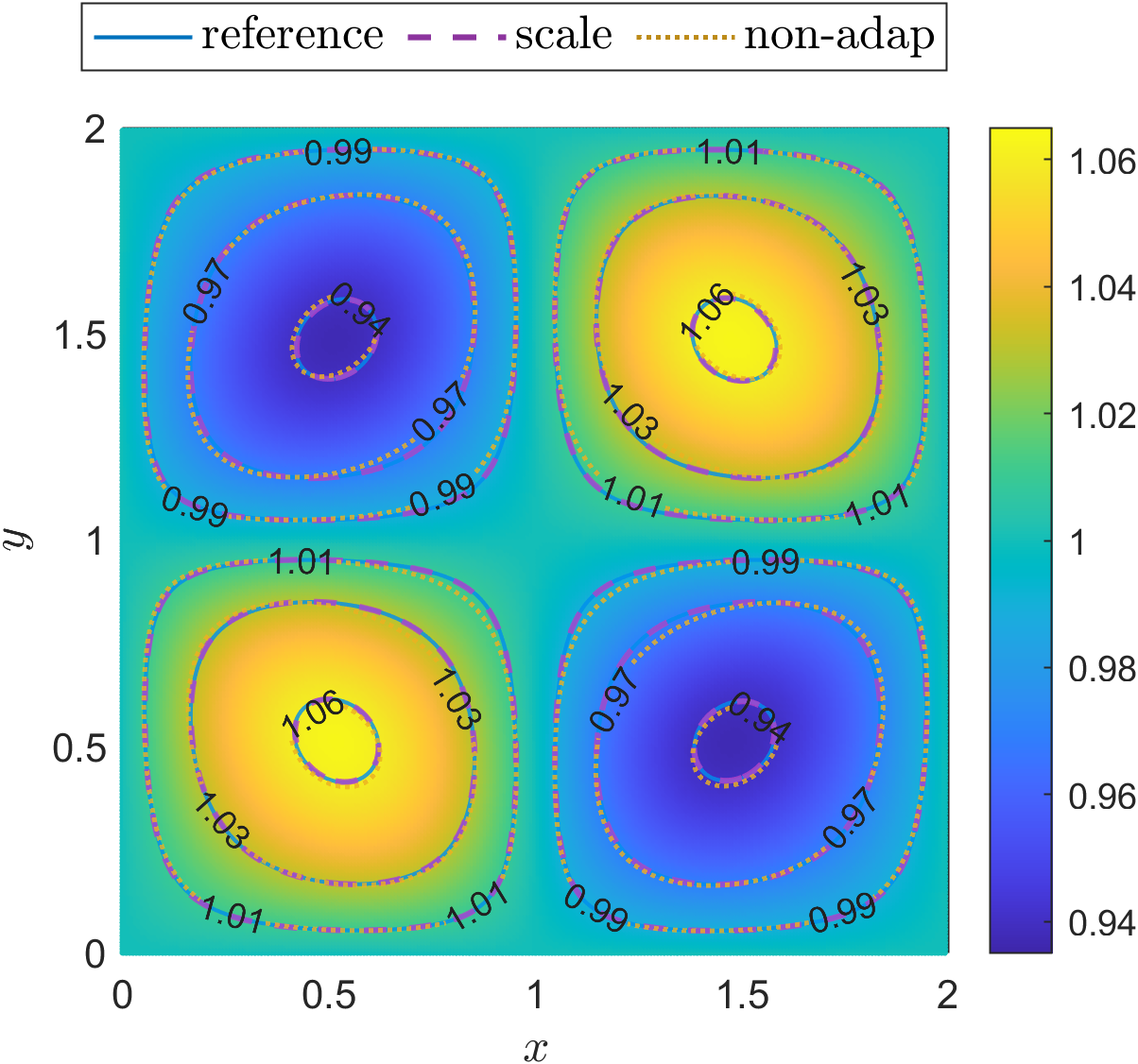}
    }\hspace{-0.7em}
    \subfigure[$\theta$, $z=0.5$\label{fig:31-1b}]{
        \includegraphics[width=0.32\linewidth]{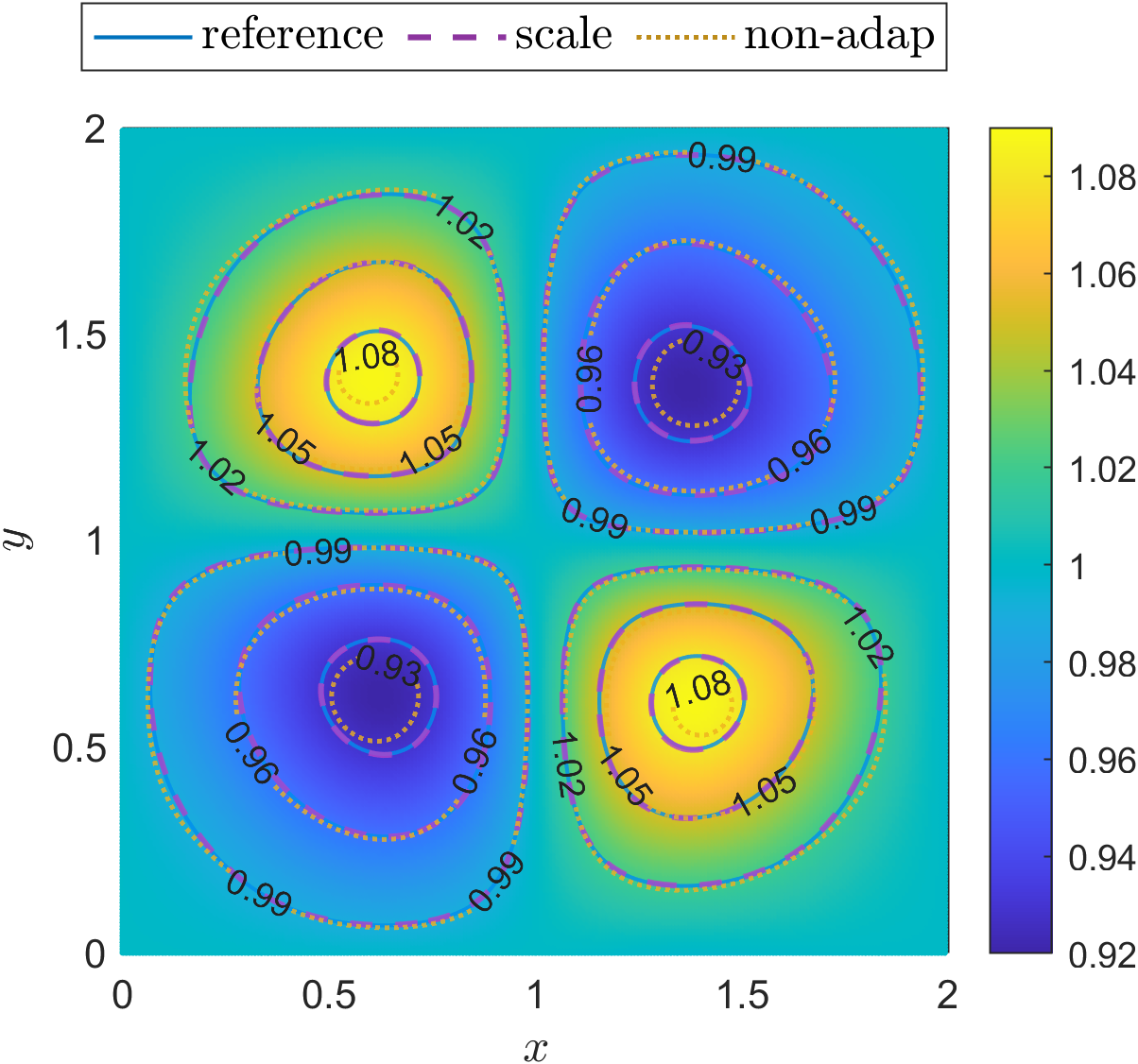}
    }\hspace{-0.7em}
    \subfigure[$q_2$, $z=0.5$\label{fig:31-1c}]{
        \includegraphics[width=0.32\linewidth]{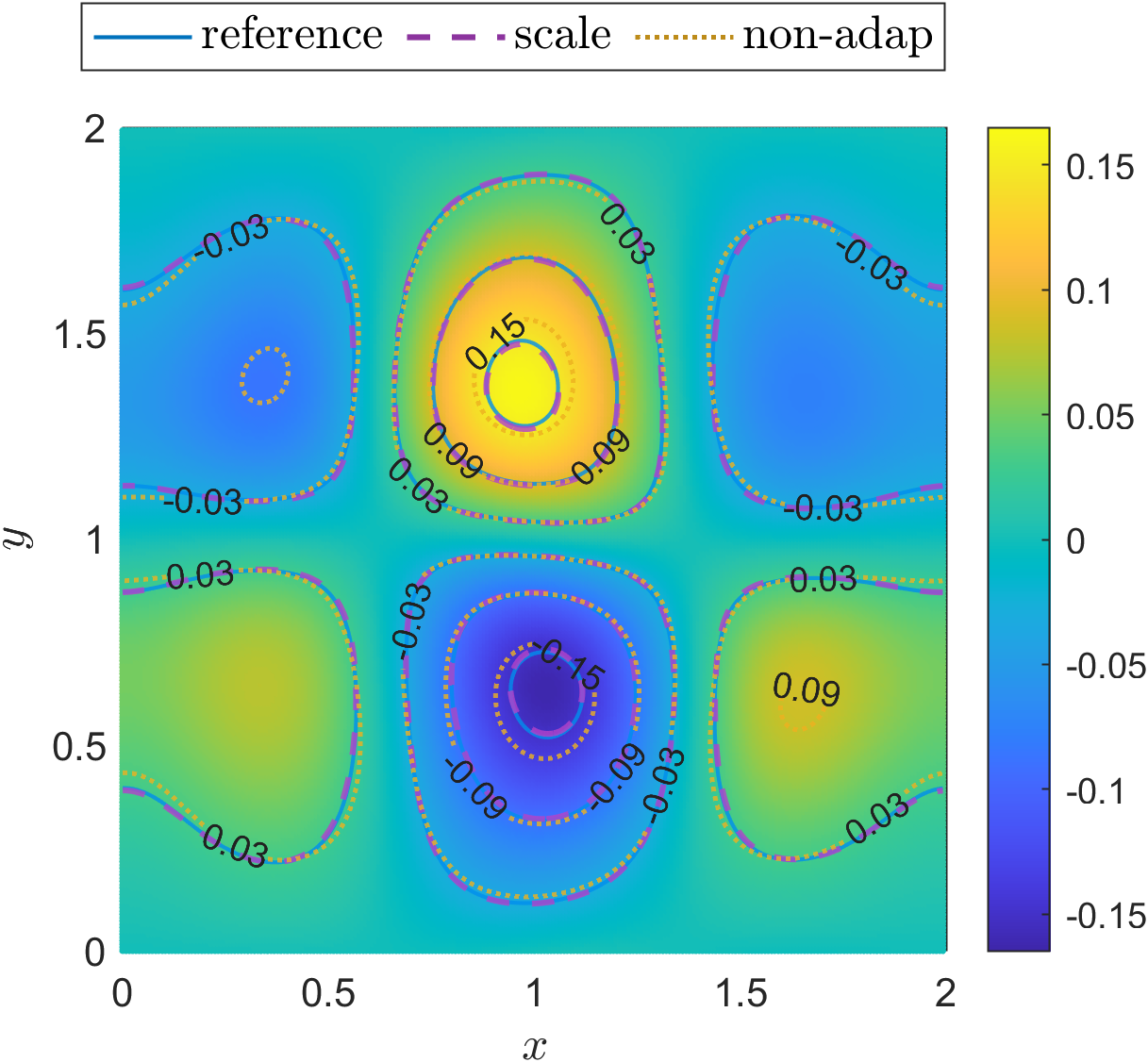}
    }
    
    \subfigure[$\rho$, $z=1.25$\label{fig:31-1d}]{
        \includegraphics[width=0.32\linewidth]{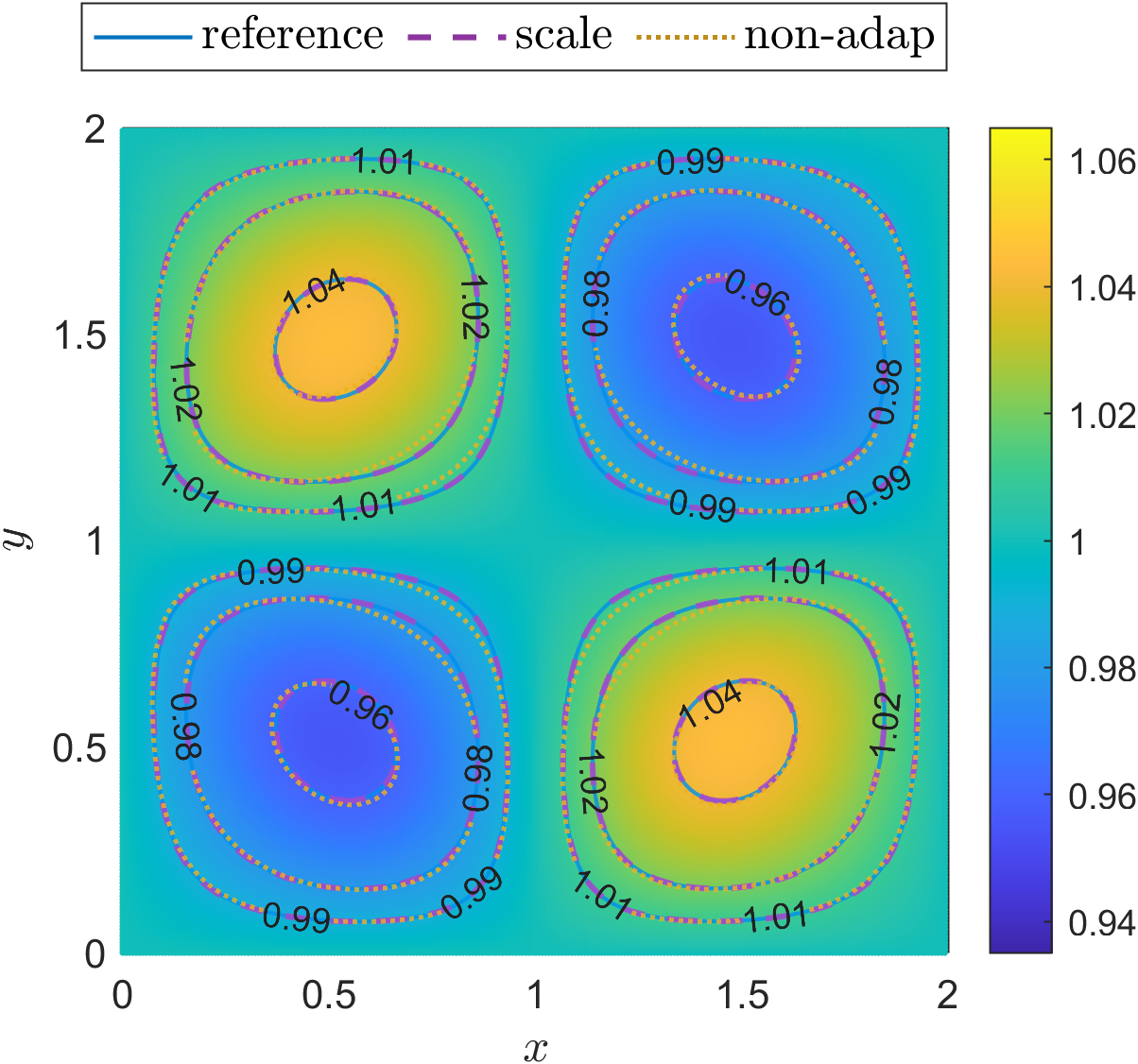}
    }\hspace{-0.7em}
    \subfigure[$\theta$, $z=1.25$\label{fig:31-1e}]{
        \includegraphics[width=0.32\linewidth]{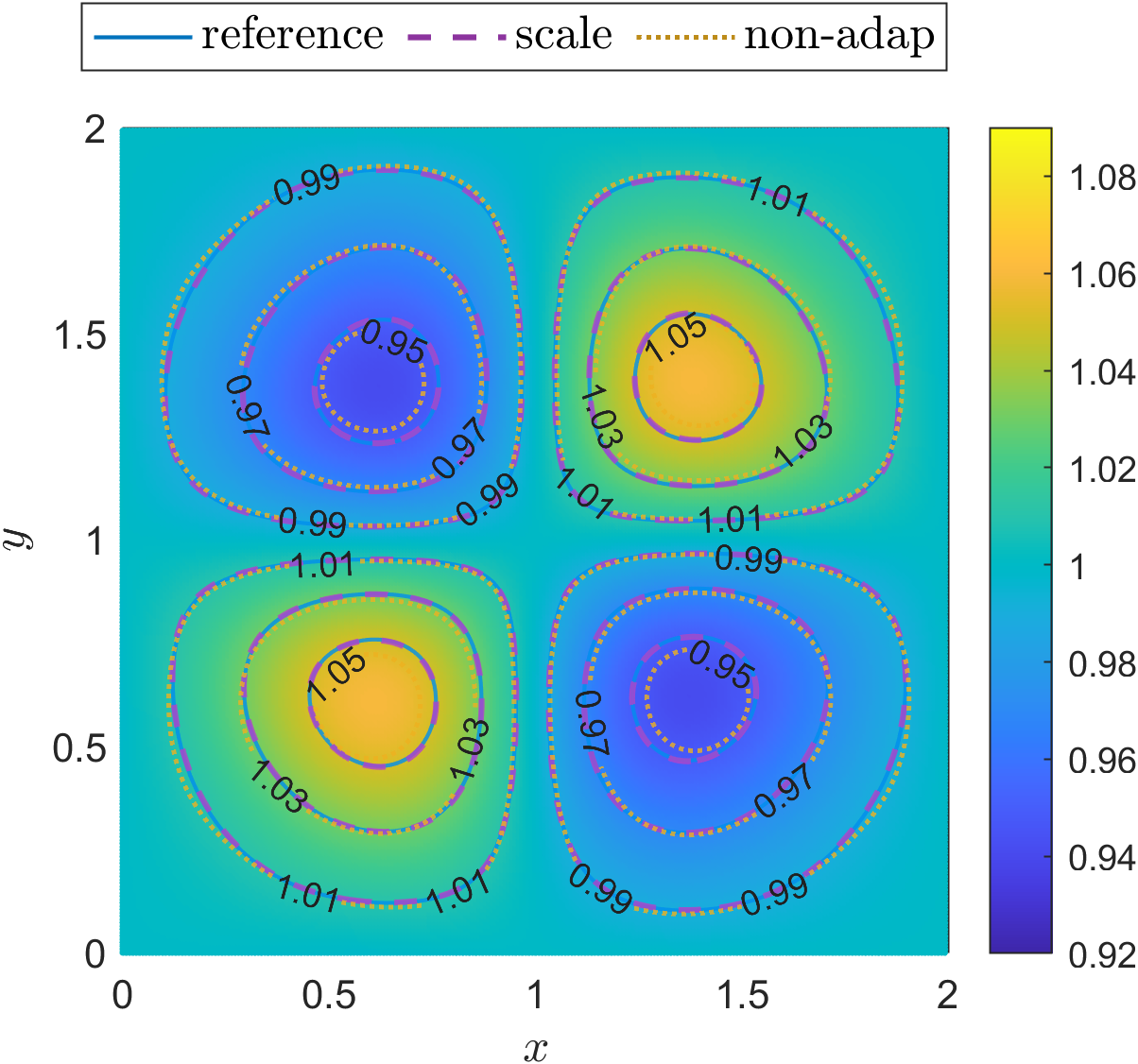}
    }\hspace{-0.7em}
    \subfigure[$q_2$, $z=1.25$\label{fig:31-1f}]{
        \includegraphics[width=0.32\linewidth]{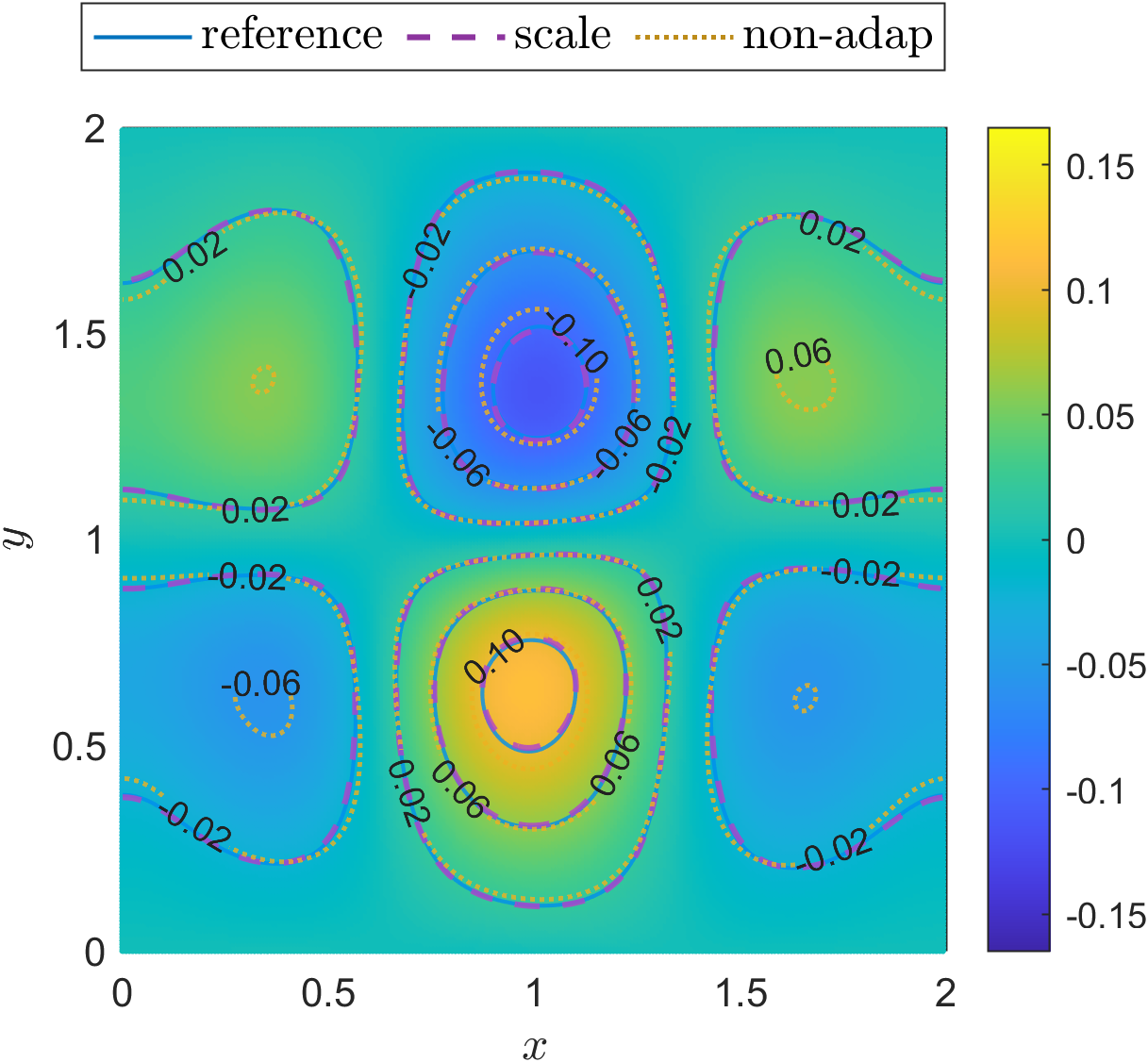}
    }
    
    \caption{(3D: Hexa-Gaussian perturbation problem in Sec.~\ref{sec:Eg31}) The numerical solution of the density $\rho$, temperature $\theta$, and heat flux $q_2$ obtained by the scaling adaptive, the non-adaptive method with $N=8$ and the reference solution at $t=0.3$. (a) Density $\rho$ at the plane $z=0.5$. (b) Temperature $\theta$ at the plane $z=0.5$. (c) Heat flux $q_2$ at the plane $z=0.5$. (d) Density $\rho$ at the plane $z=1.25$. (e) Temperature $\theta$ at the plane $z=1.25$. (f) Heat flux $q_2$ at the plane $z=1.25$.}
    \label{fig:31-1}
\end{figure}

The numerical solution of the density $\rho$, temperature $\theta$, heat flux in the $y$-axis $q_2$ on the plane $z = 0.5$, and $z = 1.25$ at $t = 0.3$ is shown in Fig.~\ref{fig:31-1}. The numerical solution by the non-adaptive method with $N = 8$, and the reference solution are both plotted in Fig.~\ref{fig:31-1}. It shows that, for all these macroscopic variables, the numerical solution by the adaptive method matches well with the reference solution, while there exist several differences between the numerical solution by the non-adaptive method and the reference solution, indicating the greater accuracy of the adaptive method. 

\begin{figure}
    \centering
    \subfigure[$L^2$ error of $\rho$ and $u_1$\label{fig:31-2a}]{
        \includegraphics[height=7.5\baselineskip]{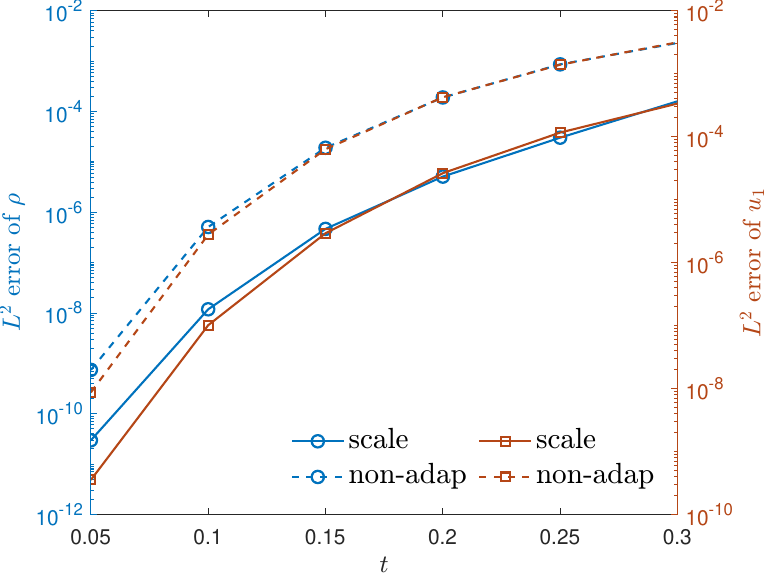}
    }\hspace{-0.7em}
    \subfigure[$L^2$ error of $\theta$ and $\sigma_{23}$\label{fig:31-2b}]{
        \includegraphics[height=7.5\baselineskip]{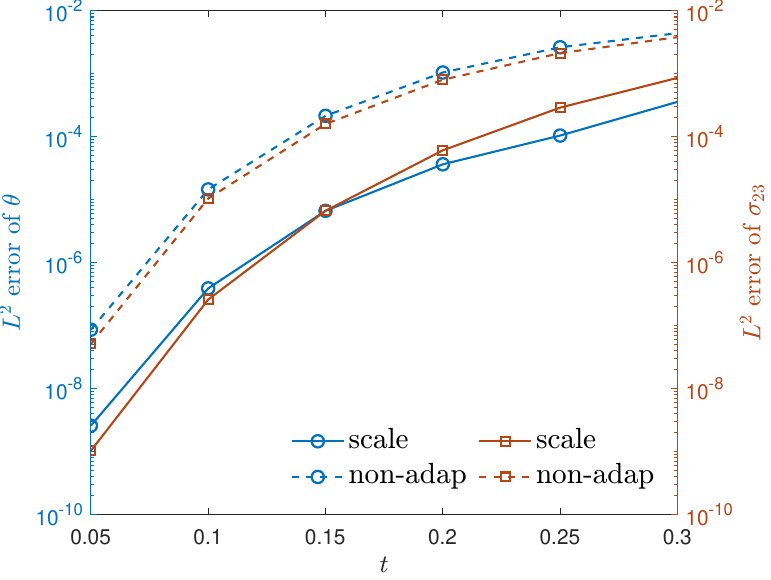}
    }\hspace{-0.7em}
    \subfigure[$L^2$ error of $q_2$ and $q_3$\label{fig:31-2c}]{
        \includegraphics[height=7.5\baselineskip]{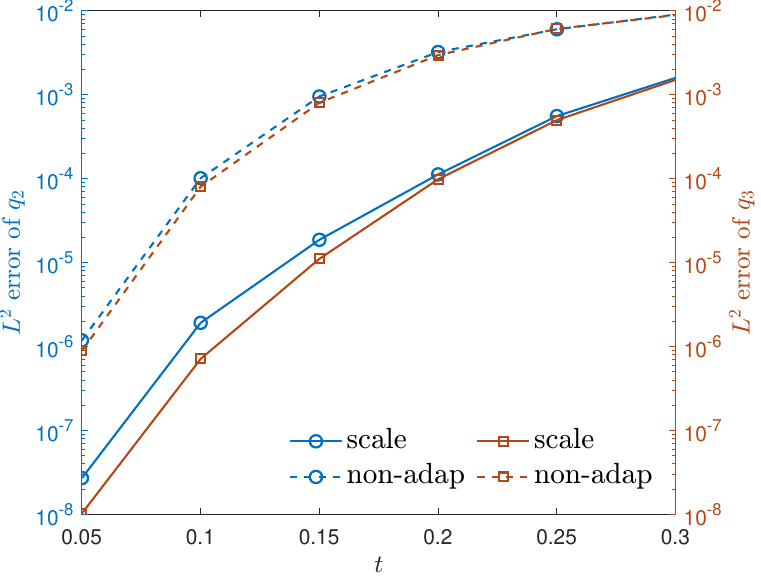}
    }
    
    \caption{(3D: Hexa-Gaussian perturbation problem in Sec.~\ref{sec:Eg31}) Comparison between the scaling adaptive and non-adaptive methods with $N=8$. (a) Evolution of the $L^2$ error between the numerical and the reference macroscopic density $\rho$ (blue) and velocity $u_1$ (red). (b) Evolution of the $L^2$ error between the numerical and the reference temperature $\theta$ (blue) and stress $\sigma_{23}$ (red). (c) Evolution of the $L^2$ error between the numerical and the reference heat flux $q_2$ (blue) and $q_3$ (red).}
    \label{fig:31-2}
\end{figure}

To show the high accuracy of the adaptive method more clearly, the evolution of the $L^2$ error for the density $\rho$, the macroscopic velocity $u_1$, temperature $\theta$, shear stress $\sigma_{23}$, the heat flux $q_2$, and $q_3$ are presented in Fig.~\ref{fig:31-2}. It reveals that for all these macroscopic variables, a significantly lower $L^2$ error is obtained for the scaling adaptive method compared to the non-adaptive method. All these demonstrate the high efficiency of the adaptive method. 

The computational costs for the scaling adaptive method and the non-adaptive method with expansion order $N=8$ are listed in Tab.~\ref{tab:11-1}. OpenMP with 32 threads is employed to achieve acceleration. The computational time for the two methods is nearly identical, which is consistent with the previous simulations.

\section{Conclusion and Outlook}\label{sec:conclusion}
In this work, an adaptive method is proposed in the framework of the Hermite spectral method for the Boltzmann equation. A $p$- and a scaling adaptive method are firstly presented for the homogeneous problems, where a frequency indicator is constructed based on the expansion coefficients for both adaptive methods. 

Subsequently, the frequency indicator is extended to the non-homogeneous problem where the Fourier spectral method is utilized for the spatial discretization. Then, a combined adaptive Hermite spectral method composed of the scaling and $p$-adaptive method is proposed to solve the Boltzmann equation. Finally, the high efficiency of this adaptive method is validated by three homogeneous and four non-homogeneous cases, respectively. For future work, we aim to extend the applications of this adaptive Hermite method to other problems with shock or Maxwell boundary conditions.


\appendix
\section{Detailed approximation of the collision operator}\label{Supp:3}
Here, a detailed version of Sec.~\ref{sec:32Q} is presented. Following the approach in \cite{Wang2019}, we first approximate the collision term using a truncated series:
\begin{align}
    &Q[f_{N}^{\bm{\zeta},\beta},f_{N}^{\bm{\zeta},\beta}](\bm{v})\approx Q_{N}^{\bm{\zeta},\beta}(\bm v) =\sum_{|k|=0}^{N}\hat{Q}_{k}^{\bm{\zeta},\beta}\mathcal{H}_{k}^{\beta}(\bm{v}-\bm{\zeta}), \\ 
    &\hat{Q}_{k}^{\bm{\zeta},\beta}=\int_{\mathbb{R}^3}H_{k}^{\beta}(\bm{v}-\bm{\zeta})Q[f_{N}^{\bm{\zeta},\beta},f_{N}^{\bm{\zeta},\beta}](\bm{v})\dd \bm{v}, \quad 0\leqslant |k|\leqslant N.
\end{align}
Substituting the expansion \eqref{eq:3fN} into the bilinear form \eqref{eq:2Q} yields a relation between the coefficients of the collision term $Q[f_{N}^{\bm{\zeta},\beta},f_{N}^{\bm{\zeta},\beta}]$ and distribution $f_{N}^{\bm{\zeta},\beta}$
\begin{equation}\label{eq:3Qhat}
    \hat{Q}_{k}^{\bm{\zeta},\beta}=\sum_{|i|=0}^{N}\sum_{|j|=0}^{N} A^{i,j;\beta}_{k}\hat{f}_{i}^{\bm{\zeta},\beta}\hat{f}_{j}^{\bm{\zeta},\beta},\qquad 0\leqslant|k|\leqslant N,
\end{equation}
where
\begin{equation}\label{eq:3A}
\begin{aligned}
    A^{i,j;\beta}_{k} = \int_{\mathbb{R}^3}\int_{\mathbb{R}^3}\int_{\mathbb{S}^2}& [\mathcal{H}_{i}^{\beta}(\bm{v}_{*}^{\prime}-\bm{\zeta})\mathcal{H}_{j}^{\beta}(\bm{v}^{\prime}-\bm{\zeta})-\mathcal{H}_{i}^{\beta}(\bm{v}_{*}-\bm{\zeta})\mathcal{H}_{j}^{\beta}(\bm{v}-\bm{\zeta})]\\
    &B(|\bm{g}|,\chi)H_{k}^{\beta}(\bm{v}-\bm{\zeta})\,\dd\bm{\omega}\dd\bm{v}_{*}\dd\bm{v}.
\end{aligned}
\end{equation}
Note that the moving center $\bm{\zeta}$ on the right-hand side of \eqref{eq:3A} can be eliminated. The ninth-order tensor $A^{i,j;\beta}_{k}$, involving an eight-dimensional integral, plays a critical role in evaluating the collision term. To avoid redundant computation, it is necessary to factor out the scaling parameter $\beta$ from \eqref{eq:3A}. Observe that
\begin{equation}
    \mathcal{H}_{k}^{\beta}(\bm{v})=\beta^{\frac{3}{2}}\mathcal{H}_{k}^{1}(\beta\bm{v}),\qquad
    H_{k}^{\beta}(\bm{v})=\beta^{\frac{3}{2}}H_{k}^{1}(\beta\bm{v}).
\end{equation}
Furthermore, the IPL model (\ref{eq:2IPL}) gives
\begin{equation}
	B\left(\left|\frac{\bm{g}}{\beta}\right|,\chi\right)=\beta^{-\frac{\eta-5}{\eta-1}}B(|\bm{g}|,\chi).
\end{equation}
After a change of variables, we obtain from \eqref{eq:3A} that
\begin{align}
    A^{i,j;\beta}_{k}&=\int_{\mathbb{R}^3}\int_{\mathbb{R}^3}\int_{\mathbb{S}^2} \beta^3 [\mathcal{H}_{i}^{1}(\bm{v}_{*}^{\prime})\mathcal{H}_{j}^{1}(\bm{v}^{\prime})-\mathcal{H}_{i}^{1}(\bm{v}_{*})\mathcal{H}_{j}^{1}(\bm{v})]\beta^{-\frac{\eta-5}{\eta-1}}B(|\bm{g}|,\chi)\,\\&\qquad\qquad\qquad\beta^{\frac{3}{2}}H_{k}^{1}(\bm{v})\,\beta^{-6}\dd\bm{\omega}\dd\bm{v}_{*}\dd\bm{v}\\&
    =\beta^{-\frac{3}{2}-\frac{\eta-5}{\eta-1}}A^{i,j;1}_{k}.
\end{align}
Using this relation, we can precompute and store the tensor $\left\lbrace A^{i,j;1}_{k}\right\rbrace$. During computation, this tensor is scaled to $\left\lbrace A^{i,j;\beta}_{k}\right\rbrace$ according to the current value of $\beta$, as adjusted by the adaptive algorithm. The details of computing the tensor $\left\lbrace A^{i,j;1}_{k}\right\rbrace$ are beyond the scope of this work. Interested readers are referred to \cite{Wang2019}.

The computational complexity of evaluating all coefficients $\hat{Q}_{k}^{\bm{\zeta},\beta}$ in \eqref{eq:3Qhat} is of $\mathcal{O}(N^9)$. In some cases, a large number of $N$ is required to achieve sufficient approximation accuracy. The resulting time and memory demands could exceed the capabilities of a single machine. To address this, the combined collision model introduced in \cite{Wang2019, Hu2020} is adopted here to help mitigate the high computational cost associated with large $N$.

Specifically, a moderate constant order $N_0$ is introduced. The core idea is to apply the full quadratic collision operator only to the coefficients $\hat{Q}_{k}^{\bm{\zeta},\beta}$ with $0\leqslant|k|\leqslant N_0$. For the remaining terms, a BGK-type operator is utilized. The resulting combined model is given by
\begin{equation}
    \hat{Q}_{k}^{\bm{\zeta},\beta}=\begin{cases}
        \displaystyle\sum_{|i|=0}^{N_0}\sum_{|j|=0}^{N_0} A^{i,j;\beta}_{k}\hat{f}_{i}^{\bm{\zeta},\beta}\hat{f}_{j}^{\bm{\zeta},\beta}, & 0\leqslant|k|\leqslant N_0,\\
        \nu_{N_0}\left(\hat{\mathcal{M}}_{k}-\hat{f}_{k}^{\bm{\zeta},\beta}\right), & N_0<|k|\leqslant N.
    \end{cases}
\end{equation}
Here, the constant $\nu_{N_0}$ denotes the spectral radius of the discrete linearized collision operator. For further details, we refer the reader to \cite{Wang2019}. The Hermite coefficients of the local Maxwellian, denoted by $\hat{\mathcal{M}}_{k}$ for $N_0<|k|\leqslant N$, can be calculated recursively in $\mathcal{O}(N^3)$ operations. To conclude, if $N\leqslant N_0$, the full collision model \eqref{eq:3Qhat} is adopted with complexity $\mathcal{O}(N^9)$. Otherwise, if $N>N_0$, the combined model \eqref{eq:5Qhat} is adopted and the complexity is reduced to a more feasible $\mathcal{O}(N_0^9+N^3)$.

Despite this simplification, accurately evaluating the collision term remains computationally expensive. However, it is known that by appropriately selecting the Hermite basis functions, one can achieve equal or better accuracy with a smaller expansion order $N$. We emphasize that a significant reduction in computation time can be realized when $N<N_0$. This is precisely the goal of the adaptive technique introduced in the next section.

\section{Proof of Proposition \ref{prop:1}}\label{Supp:1}
Proofs of Eqs. (\ref{eq:prop1}--\ref{eq:prop4}) are as follows.

\begin{itemize}[leftmargin=4em]
    \item \eqref{eq:prop1} can be directly deduced from the orthogonality and the odd-even parity of $H_l^{\beta'}$ and $\mathcal{H}_{k}^{\beta}$.
    
    \item Notice that
    \begin{equation*}
        \frac{\dd \mathcal{H}_k^{\beta}(v)}{\dd v}=-\beta\sqrt{k+1}\mathcal{H}_{k+1}^{\beta}(v),\quad \frac{\dd H_{l+1}^{\beta}(v)}{\dd v}=\beta\sqrt{l+1}H_{l}^{\beta}(v),\quad k,l\geqslant0.
    \end{equation*}
    The first equation of \eqref{eq:prop2} is derived using integration by parts
    \begin{align*}
        T_{l,k}&=\int_{\mathbb{R}} \mathcal{H}_k^{\beta}(v)H_l^{\beta^\prime}(v)\dd v\\
        &=\left.\frac{1}{\beta^{\prime}\sqrt{l+1}}\mathcal{H}_k^{\beta}(v)H_{l+1}^{\beta^\prime}(v)\right|_{-\infty}^{+\infty}-\int_{\mathbb{R}}\frac{1}{\beta^{\prime}\sqrt{l+1}}H_{l+1}^{\beta^\prime}(v)\dd\mathcal{H}_k^{\beta}(v)\\
        &=\frac{\beta}{\beta^\prime}\sqrt{\frac{k+1}{l+1}}\int_{\mathbb{R}} H_{l+1}^{\beta^\prime}(v)\mathcal{H}_{k+1}^{\beta}(v)\dd v\\
        &=\frac{\beta}{\beta^{\prime}}\sqrt{\frac{k+1}{l+1}}T_{l+1,k+1}.
    \end{align*}
    For the latter, the initial term $T_{0,0}$ can be checked by substituting the explicit form of $\mathcal{H}_0^{\beta}(v)$ and $H_0^{\beta^\prime}(v)$
    \begin{equation*}
        T_{0,0}=\int_{\mathbb{R}}\sqrt{\frac{\beta\beta'}{2\pi}}\ee^{-\frac{\beta^2 v^2}{2}}\dd v=\sqrt{\frac{\beta'}{\beta}},
    \end{equation*}
    Then the cases of $l>0$ can be obtained using the first equation recursively.
    
    \item Due to the three-term recurrence relation of the Hermite polynomials 
    \begin{align*}
        H_{l+1}^{\beta}(v)=\frac{\beta v}{\sqrt{l+1}}H_{l}^{\beta}(v)-\sqrt{\frac{l}{l+1}}H_{l-1}^{\beta}(v),
    \end{align*}
    it holds that 
    \begin{align*}
        T_{l,k}&=\int\mathcal{H}_k^{\beta}(v)H_l^{\beta^{\prime}}(v)\dd v \\
        &=\int\mathcal{H}_k^{\beta}(v)\left[ \frac{\beta^{\prime}v}{\sqrt{l}}H_{l-1}^{\beta^{\prime}}(v)-\sqrt{\frac{l-1}{l}}H_{l-2}^{\beta^{\prime}}(v)\right] \dd v\\
        &=\frac{\beta^{\prime}}{\sqrt{l}}\int v\mathcal{H}_k^{\beta}(v)H_{l-1}^{\beta^{\prime}}(v)\dd v-\sqrt{\frac{l-1}{l}}T_{l-2,k}\\
        &=\frac{\beta^{\prime}}{\beta}\frac{1}{\sqrt{l}}\int\beta v\mathcal{H}_k^{\beta}(v)H_{l-1}^{\beta^{\prime}}(v)\dd v-\sqrt{\frac{l-1}{l}}T_{l-2,k}\\
        &=\frac{\beta^{\prime}}{\beta}\frac{1}{\sqrt{l}}\int\left[\sqrt{k+1}\mathcal{H}_{k+1}^{\beta}(v)+\sqrt{k}\mathcal{H}_{k-1}^{\beta}(v)\right] H_{l-1}^{\beta^{\prime}}(v)\dd v-\sqrt{\frac{l-1}{l}}T_{l-2,k}\\
        &=\frac{\beta^{\prime}}{\beta}\sqrt{\frac{k+1}{l}}T_{l-1,k+1}+\frac{\beta^{\prime}}{\beta}\sqrt{\frac{k}{l}}T_{l-1,k-1}-\sqrt{\frac{l-1}{l}}T_{l-2,k}.
    \end{align*}
    \item Combining \eqref{eq:prop2} and \eqref{eq:prop3}, we have
    \begin{equation*}
        T_{l,k}=\frac{\beta^{\prime}}{\beta}\sqrt{\frac{k+1}{l}}T_{l-1,k+1}+\frac{k}{l}T_{l,k}-\frac{\beta}{\beta^{\prime}}\sqrt{\frac{k+1}{l}}T_{l-1,k+1}.
    \end{equation*}
    Simplifying this gives \eqref{eq:prop4}.
\end{itemize}
So far, we have completed the proof of Prop.~\ref{prop:1}.

\section{Derivation of advection term}\label{Supp:2}
To compute the advection term, we need the following relation
\begin{equation*}
    v_1 \mathcal{H}_{k}^{\beta}(\bm{v}-\bm{\zeta}) =\frac{\sqrt{k_1+1}}{\beta}\mathcal{H}_{k+e_1}^{\beta}(\bm{v}-\bm{\zeta}) +\zeta_1 \mathcal{H}_{k}^{\beta}(\bm{v}-\bm{\zeta}) +\frac{\sqrt{k_1}}{\beta}\mathcal{H}_{k-e_1}^{\beta}(\bm{v}-\bm{\zeta}).
\end{equation*}
Then
\begin{align*}
    &v_1 f_{N,M}^{\bzeta,\beta}(t,x_j,\bm v)=\sum_{|k|=0}^{N+1}\hat{F}_{k}(t,x_j) \mathcal{H}_{k}^{\beta}(\bm{v}-\bm{\zeta}),\\
    &\hat{F}_{k}=\frac{\sqrt{k_1}}{\beta}\hat{f}_{k-e_1}^{\bzeta,\beta}+\zeta_1 \hat{f}_{k}^{\bzeta,\beta}+\frac{\sqrt{k_1+1}}{\beta}\hat{f}_{k+e_1}^{\bzeta,\beta},\quad 0\leqslant|k|\leqslant N+1,
\end{align*}
with $\hat{f}_{k}^{\bzeta,\beta}\coloneq0$ for $|k|\notin[0,N]$. Note that
\begin{equation*}
    \hat{F}_{k}(t,\cdot)\in\Span\left\lbrace E_{l}\mid -M/2\leqslant l\leqslant M/2\right\rbrace.
\end{equation*}
The derivatives can be evaluated exactly in the frequency domain. That is
\begin{align*}
    &\hat{F}_{k}(t,x_j)=\sum_{l=-M/2}^{M/2}\hat{G}_{k,l}(t)E_{l}(x_j),\\ 
    &\hat{G}_{k,l}(t)=\frac{\sqrt{k_1}}{\beta}\hat{g}_{k-e_1,l}^{\bzeta,\beta}(t)+\zeta_1 \hat{g}_{k,l}^{\bzeta,\beta}(t)+\frac{\sqrt{k_1+1}}{\beta}\hat{g}_{k+e_d,l}^{\bzeta,\beta}(t).
\end{align*}
Then we have
\begin{equation*}
    \frac{\pp\hat{F}_{k}}{\pp x}(t,x_j)=\sum_{l=-M/2}^{M/2}\frac{2\pi\ii l}{L}\hat{G}_{k,l}(t)E_{l}(x_j).
\end{equation*}
which is the same as \eqref{eq:5adv} and \eqref{eq:5h_hat}.

\bibliographystyle{siamplain}
\bibliography{references}

\begin{thebibliography}{10}

\bibitem{BenAbdallah1996}
{\sc N.~Abdallah, P.~Degond, and S.~Genieys}, {\em An energy-transport model
  for semiconductors derived from the {Boltzmann} equation}, J. Stat. Phys., 84
  (1996), pp.~205--231.

\bibitem{Bird1994}
{\sc G.~Bird}, {\em Molecular Gas Dynamics and the Direct Simulation of Gas
  Flows}, Clarendon Press, Oxford, 1994.

\bibitem{Bobylev1975}
{\sc A.~Bobylev}, {\em Exact solutions of the {Boltzmann} equation}, Dokl.
  Akad. Nauk. SSSR, 225 (1975), pp.~1296--1299.
\newblock In Russian.

\bibitem{Cai2021}
{\sc Z.~Cai}, {\em Moment method as a numerical solver: Challenge from shock
  structure problems}, J. Comput. Phys., 444 (2021), p.~110593.

\bibitem{Cai2014d}
{\sc Z.~Cai, Y.~Fan, and R.~Li}, {\em Globally hyperbolic regularization of
  {Grad}’s moment system}, Commun. Pure Appl. Math., 67 (2014), pp.~464--518.

\bibitem{Chorin1968}
{\sc A.~Chorin}, {\em Numerical solution of the {Navier-Stokes} equations},
  Math. Comput., 22 (1968), pp.~745--762.

\bibitem{Chou2023}
{\sc T.~Chou, S.~Shao, and M.~Xia}, {\em Adaptive {Hermite} spectral methods in
  unbounded domains}, Appl. Numer. Math., 183 (2023), pp.~201--220.

\bibitem{Deng2025}
{\sc Y.~Deng, S.~Shao, A.~Mogilner, and M.~Xia}, {\em Adaptive
  hyperbolic-cross-space mapped {Jacobi} method on unbounded domains with
  applications to solving multidimensional spatiotemporal integrodifferential
  equations}, J. Comput. Phys., 520 (2025), p.~113492.

\bibitem{Dimarco2018}
{\sc G.~Dimarco, R.~Loubère, J.~Narski, and T.~Rey}, {\em An efficient
  numerical method for solving the {Boltzmann} equation in multidimensions}, J.
  Comput. Phys., 353 (2018), pp.~46--81.

\bibitem{Dimarco2014}
{\sc G.~Dimarco and L.~Pareschi}, {\em Numerical methods for kinetic
  equations}, Acta Numer., 23 (2014), pp.~369--520.

\bibitem{Filbet2022}
{\sc F.~Filbet and T.~Xiong}, {\em Conservative discontinuous
  {Galerkin/Hermite} spectral method for the {Vlasov-Poisson} system}, Commun.
  Appl. Math. Comput., 4 (2022), pp.~34--59.

\bibitem{Goldstein1989}
{\sc D.~Goldstein, B.~Sturtevant, and J.~Broadwell}, {\em Investigations of the
  motion of discrete-velocity gases}, Prog. Astronaut. Aeronaut., 117 (1989),
  pp.~100--117.

\bibitem{Hu2024a}
{\sc H.~Hu and H.~Yu}, {\em Scaling optimized {Hermite} approximation methods},
  arXiv:2412.08044,  (2024).

\bibitem{Hu2022a}
{\sc J.~Hu, X.~Huang, J.~Shen, and H.~Yang}, {\em A fast {Petrov-Galerkin}
  spectral method for the multidimensional {Boltzmann} equation using mapped
  {Chebyshev} functions}, SIAM J. Sci. Comput., 44 (2022), pp.~A1497--A1524.

\bibitem{Hu2020}
{\sc Z.~Hu, Z.~Cai, and Y.~Wang}, {\em Numerical simulation of microflows using
  {Hermite} spectral methods}, SIAM J. Sci. Comput., 42 (2020), pp.~B105--B134.

\bibitem{Issan2024}
{\sc O.~Issan, O.~Koshkarov, F.~Halpern, B.~Kramer, and G.~Delzanno}, {\em
  Anti-symmetric and positivity preserving formulation of a spectral method for
  {Vlasov-Poisson} equations}, J. Comput. Phys., 514 (2024), p.~113263.

\bibitem{Krook1977}
{\sc M.~Krook and T.~Wu}, {\em Exact solutions of the {Boltzmann} equation},
  Phys. Fluids, 20 (1977), pp.~1589--1595.

\bibitem{Pagliantini2023}
{\sc C.~Pagliantini, G.~Delzanno, and S.~Markidis}, {\em Physics-based
  adaptivity of a spectral method for the {Vlasov-Poisson} equations based on
  the asymmetrically-weighted {Hermite} expansion in velocity space}, J.
  Comput. Phys., 488 (2023), p.~112252.

\bibitem{Pareschi1996}
{\sc L.~Pareschi and B.~Perthame}, {\em A {Fourier} spectral method for
  homogeneous {Boltzmann} equations}, Transport Theory Statist. Phys., 25
  (1996), pp.~369--382.

\bibitem{Pareschi2000b}
{\sc L.~Pareschi and G.~Russo}, {\em Numerical solution of the {Boltzmann}
  equation {I}: Spectrally accurate approximation of the collision operator},
  SIAM J. Numer. Anal., 37 (2000), pp.~1217--1245.

\bibitem{Shen2011}
{\sc J.~Shen, T.~Tang, and L.-L. Wang}, {\em Spectral Methods: Algorithms,
  Analysis and Applications}, Springer, Berlin, Heidelberg, 2011.

\bibitem{Shen2009}
{\sc J.~Shen and L.-L. Wang}, {\em Some recent advances on spectral methods for
  unbounded domains}, Commun. Comput. Phys., 5 (2009), pp.~195--241.

\bibitem{Shu1988}
{\sc C.-W. Shu and S.~Osher}, {\em Efficient implementation of essentially
  non-oscillatory shock-capturing schemes}, J. Comput. Phys., 77 (1988),
  pp.~439--471.

\bibitem{Tang1993}
{\sc T.~Tang}, {\em The {Hermite} spectral method for {Gaussian}-type
  functions}, SIAM J. Sci. Comput., 14 (1993), pp.~594--606.

\bibitem{Vencels2015}
{\sc J.~Vencels, G.~Delzanno, A.~Johnson, I.~Peng, E.~Laure, and S.~Markidis},
  {\em Spectral solver for multi-scale plasma physics simulations with
  dynamically adaptive number of moments}, Procedia Comput. Sci., 51 (2015),
  pp.~1148--1157.

\bibitem{Wang2019}
{\sc Y.~Wang and Z.~Cai}, {\em Approximation of the {Boltzmann} collision
  operator based on {Hermite} spectral method}, J. Comput. Phys., 397 (2019),
  p.~108815.

\bibitem{Xia2021}
{\sc M.~Xia, S.~Shao, and T.~Chou}, {\em Efficient scaling and moving
  techniques for spectral methods in unbounded domains}, SIAM J. Sci. Comput.,
  43 (2021), pp.~A3244--A3268.

\bibitem{Xia2021a}
{\sc M.~Xia, S.~Shao, and T.~Chou}, {\em A frequency-dependent $p$-adaptive
  technique for spectral methods}, J. Comput. Phys., 446 (2021), p.~110627.

\bibitem{Xiao2020}
{\sc T.~Xiao, C.~Liu, K.~Xu, and Q.~Cai}, {\em A velocity-space adaptive
  unified gas kinetic scheme for continuum and rarefied flows}, J. Comput.
  Phys., 415 (2020), p.~109535.

\end{thebibliography}

\end{document}